\documentclass[a4paper,12pt]{article}%
\usepackage{amssymb}
\usepackage{amsmath}
\usepackage{latexsym}
\usepackage{graphicx}
\usepackage[latin1]{inputenc}
\usepackage{amsfonts}
\usepackage{hyperref}
\usepackage{amsthm}
\usepackage{graphicx}
\usepackage{color}
\usepackage{verbatim}%
\providecommand{\U}[1]{\protect\rule{.1in}{.1in}}
\font\teneusb=eusb10 \font\seveneusb=eusb7 \font\fiveeusb=eusb5
\newfam\eusbfam \textfont\eusbfam=\teneusb
\scriptfont\eusbfam=\seveneusb \scriptscriptfont\eusbfam=\fiveeusb

\font\tenbifull=cmmib10
\font\tenbimed=cmmib7
\font\tenbismall=cmmib5
\mathchardef\bbGamma="7000 \mathchardef\bbDelta="7001
\mathchardef\bbPhi="7002 \mathchardef\bbAlpha="7003
\mathchardef\bbXi="7004 \mathchardef\bbPi="7005
\mathchardef\bbSigma="7006 \mathchardef\bbUpsilon="7007
\mathchardef\bbTheta="7008 \mathchardef\bbPsi="7009
\mathchardef\bbOmega="700A \mathchardef\bbalpha="710B
\mathchardef\bbbeta="710C \mathchardef\bbgamma="710D
\mathchardef\bbdelta="710E \mathchardef\bbepsilon="710F
\mathchardef\bbzeta="7110 \mathchardef\bbeta="7111
\mathchardef\bbtheta="7112 \mathchardef\bbiota="7113
\mathchardef\bbkappa="7114 \mathchardef\bblambda="7115
\mathchardef\bbmu="7116 \mathchardef\bbnu="7117
\mathchardef\bbxi="7118 \mathchardef\bbpi="7119
\mathchardef\bbrho="711A \mathchardef\bbsigma="711B
\mathchardef\bbtau="711C \mathchardef\bbupsilon="711D
\mathchardef\bbphi="711E \mathchardef\bbchi="711F
\mathchardef\bbpsi="7120 \mathchardef\bbomega="7121
\mathchardef\bbvarepsilon="7122 \mathchardef\bbvartheta="7123
\mathchardef\bbvarpi="7124 \mathchardef\bbvarrho="7125
\mathchardef\bbvarsigma="7126 \mathchardef\bbvarphi="7127

\newcommand{\N}{{\rm I}\kern-0.18em{\rm N}}

\newcommand{\h}{{\rm I}\kern-0.18em{\rm H}}
\newcommand{\K}{{\rm I}\kern-0.18em{\rm K}}

\newcommand{\Z}{{\rm Z}\kern-0.34em{\rm Z}}

\newcommand{\1}{{\rm 1}\kern-0.22em{\rm I}}

\newtheorem{thm}{Theorem}[section]

\newtheorem{ex}{Example}[section]

\newtheorem{prop}{Proposition}[section]
\newtheorem{cor}{Corollary}[section]

\newtheorem{lem}{Lemma}[section]
\newtheorem{rem}{Remark}[section]

\numberwithin{equation}{section}

\newcounter{eqroman}
\begin{document}

\title{A conditional limit theorem for random walks under extreme deviation}
\author{Michel Broniatowski$^{(1)}$ and Zhangsheng Cao$^{(1)}$\\LSTA, Université Paris 6}
\maketitle

\begin{abstract}
This paper explores a conditional \ Gibbs theorem for a random walkinduced by
i.i.d. $(X_{1},..,X_{n})$ conditioned on an extreme deviation of its sum
$(S_{1}^{n}=na_{n})$ or $(S_{1}^{n}>na_{n})$ where $a_{n}\rightarrow\infty$.
It is proved that when the summands have light tails with some additional
regulatity property, then the asymptotic conditional distribution of $X_{1}$
can be approximated in variation norm by the tilted distribution at point
$a_{n}$ , extending therefore the classical LDP case.

\end{abstract}

\section{\bigskip Introduction}

Let $X_{1}^{n}:=\left(  X_{1},..,X_{n}\right)  $ denote $n$ independent
unbounded real valued random variables and $S_{1}^{n}:=X_{1}+..+X_{n}$ denote
their sum. The purpose of this paper is to explore the limit distribution of
the generic variable $X_{1}$ conditioned on extreme deviations (ED) pertaining
to $S_{1}^{n}.$ By extreme deviation we mean that $S_{1}^{n}/n$ is supposed to
take values which are going to infinity as $n$ increases. Obviously such
events are of infinitesimal probability. Our interest in this question stems
from a first result which assesses that under appropriate conditions, when the
sequence $a_{n}$ is such that
\[
\lim_{n\rightarrow\infty}a_{n}=\infty
\]
then there exists a sequence $\varepsilon_{n}$ which tends to $0$ as $n$ tends
to infinity such that
\begin{equation}
\lim_{n\rightarrow\infty}P\left(  \left.  \cap_{i=1}^{n}\left(  X_{i}%
\in\left(  a_{n}-\varepsilon_{n},a_{n}+\varepsilon_{n}\right)  \right)
\right\vert S_{1}^{n}/n>a_{n}\right)  =1 \label{democracy}%
\end{equation}
which is to say that when the empirical mean takes exceedingly large values,
then all the summands share the same behaviour. This result obviously requires
a number of hypotheses, which we simply quote as \lq\lq light tails" type. We
refer to \cite{BoniaCao} for this result and the connection with earlier
related works.

The above result is clearly to be put in relation with the so-called Gibbs
conditional Principle which we recall briefly in its simplest form.

Consider the case when the sequence $a_{n}$ $=a$ is constant with value larger
than the expectation of $X_{1}.$ Hence we consider the behaviour of the
summands when $\left(  S_{1}^{n}/n>a\right)  $ , under a large deviation (LD)
condition about the empirical mean. The asymptotic conditional distribution of
$X_{1}$ given $\left(  S_{1}^{n}/n>a\right)  $ is the well known tilted
distribution of $P_{X}$ with parameter $t$ associated to $a.$ Let us introduce
some notation to put this in light. The hypotheses to be stated now together
with notation are kept throughout the entire paper.

It will be assumed that $P_{X}$ , which is the distribution of $X_{1}$, has a
density $p$ with respect to the Lebesgue measure on $\mathbb{R}$. The fact
that $X_{1}$ has a light tail is captured in the hypothesis that $X_{1}$ has a
moment generating function
\[
\Phi(t):=E\exp tX_{1}%
\]
which is finite in a non void neighborhood $\mathcal{N}$ of $0.$ This fact is
usually refered to as a Cramer type condition.

Defined on $\mathcal{N}$ are the following functions. The functions
\[
t\rightarrow m(t):=\frac{d}{dt}\log\Phi(t)
\]%
\[
t\rightarrow s^{2}(t):=\frac{d}{dt}m(t)
\]

\[
t\rightarrow\mu_{j}(t):=\frac{d}{dt}s^{2}(t)\text{ \ , \ }j=3,4
\]
are the expectation and the three first centered moments of the r.v.
$\mathcal{X}_{t}$ with density
\[
\pi_{t}(x):=\frac{\exp tx}{\Phi(t)}p(x)
\]
which is defined on $\mathbb{R}$ and which is the tilted density with
parameter $t.$ When $\Phi$ is steep, meaning that
\[
\lim_{t\rightarrow t^{+}}m(t)=\infty
\]
where $t^{+}:=ess\sup\mathcal{N}$ then $m$ parametrizes the convex hull of the
support of $P_{X}.$ We refer to Barndorff-Nielsen (1978) for those properties.
As a consequence of this fact, for all $a$ in the support of $P_{X}$, it will
be convenient to define
\[
\pi^{a}=\pi_{t}%
\]
where $a$ is the unique solution of the equation $m(t)=a.$

We now come to some remark on the Gibbs conditional principle in the standard
above setting. A phrasing of this principle is:

As $n$ tends to infinity the conditional distribution of $X_{1}$ given
$\left(  S_{1}^{n}/n>a\right)  $ is $\Pi^{a},$ the distribution with density
$\pi^{a}.$

Indeed we prefer to state Gibbs principle in a form where the conditioning
event is a point condition $\left(  S_{1}^{n}/n=a\right)  .$ The conditional
distribution of $X_{1}$ given $\left(  S_{1}^{n}/n=a\right)  $ is a well
defined distribution and Gibbs conditional principle states that this
conditional distribution converges to $\Pi^{a}$ as $n$ tends to infinity. In
both settings, this convergence holds in total variation norm. We refer to
\cite{Diaconis1} for the local form of the conditioning event; we will mostly
be interested in the extension of this form in the present paper.

For all $\alpha$ (depending on $n$ or not) we will denote $p_{\alpha}$ \ the
density of the random vector $X_{1}^{k}$ conditioned upon the local event
$\left(  S_{1}^{n}=n\alpha\right)  .$ The notation $p_{\alpha}\left(
X_{1}^{k}=x_{1}^{k}\right)  $ is sometimes used to \ denote the value of the
density $p_{\alpha}$ at point $x_{1}^{k}.$ The same notation is used xhen
$X_{1},..,X_{n}$ are sampled under some $\Pi^{\alpha}$, namely $\pi^{\alpha
}(X_{1}^{k}=x_{1}^{k}).$

In \cite{Bronia} some extension of the above Gibbs principle has been
obtained. When $a_{n}$ $=a>EX_{1}$ a second order term provides a sharpening
of the conditioned Gibbs principle, stating that
\begin{equation}
\lim_{n\rightarrow\infty}\int\left\vert p_{a}(x)-g_{a}(x)\right\vert dx)=0
\label{BRCA}%
\end{equation}
where

\bigskip%

\begin{equation}
g_{a}(x):=Cp(x)\mathfrak{n}\left(  a,s_{n}^{2},x\right)  . \label{Ga}%
\end{equation}
Hereabove $\mathfrak{n}\left(  a,s_{n},x\right)  $ denotes the normal density
function at point $x$ with expectation $a$, with variance $s_{n}^{2}$, and
$s_{n}^{2}:=s^{2}(t)(n-1).$ In the above display, $C$ is a normalizing
constant. Obviously developing in this display yields
\[
g_{a}(x)=\pi^{a}(x)\left(  1+o(1)\right)
\]
which proves that (\ref{BRCA}) is a weak form of Gibbs principle, with some
improvement due to the second order term.

The paper is organized as follows. Notation and hypotheses are stated in
Section 2 , along with some necessary facts from asymptotic analysis in the
context of light tailed densities. Section 3 provides a local Gibbs
conditional principle under EDP, namely producing the approximation of the
conditional density of $X_{1},..,X_{k}$ conditionally on $\left(  \left(
1/n\right)  \left(  X_{1}+..+X_{n}\right)  =a_{n}\right)  $ for sequences
$a_{n}$ which tend to infinity, and where $k$ is fixed, independent on $n.$
The approximation is local. This result is extended in Section 4 to typical
paths under the conditional sampling scheme, which in turn provides the
approximation in variation norm for the conditional distribution; in this
extension, $k$ is equal 1, although the result clearly also holds for fixed
$k>1.$ The method used here follows closely the approach by \cite{Bronia}.
Discussion of the differences between the Gibbs principles in LDP and EDP are
discussed. Section 5 states similar results in the case when the conditioning
event is $\left(  \left(  1/n\right)  \left(  X_{1}+..+X_{n}\right)
>a_{n}\right)  $.

The main tools to be used come from asymptotic analysis and local limit
theorems, developped from \cite{Feller} and \cite{Bingham}; we also have
borrowed a number of arguments from \cite{Nagaev}. A number of technical
lemmas have been postponed to the appendix. \bigskip

\section{Notation and hypotheses}

In this paper, we consider the uniformly bounded density function $p(x)$
\begin{equation}
p(x)=c\exp\Big(-\big(g(x)-q(x)\big)\Big)\qquad x\in\mathbb{R}_{+},
\label{densityFunction}%
\end{equation}
where $c$ is some positive normalized constant. Define $h(x):=g^{\prime}(x) $.
We assume that for some And there exists some positive constant $\vartheta$ ,
for large $x$, it holds
\begin{equation}
\sup_{|v-x|<\vartheta x}|q(v)|\leq\frac{1}{\sqrt{xh(x)}}.
\label{densityFunction01}%
\end{equation}
The function $g$ is positive and satisfies
\begin{equation}
\frac{g(x)}{x}\longrightarrow\infty,\qquad x\rightarrow\infty.
\label{3section101}%
\end{equation}

Not all positive $g$'s satisfying $(\ref{3section101})$ are adapted to our
purpose. Regular functions $g$ are defined as follows. We define firstly a
subclass $R_{0}$ of the family of \emph{slowly varying} function. A function
$l$ belongs to $R_{0}$if it can be represented as
\begin{equation}
l(x)=\exp\Big(\int_{1}^{x}\frac{\epsilon(u)}{u}du\Big),\qquad x\geq1,
\label{3section102}%
\end{equation}
where $\epsilon(x)$ is twice differentiable and $\epsilon(x)\rightarrow0$ as
$x\rightarrow\infty$.

We follow the line of Juszczak and Nagaev $\cite{Nagaev}$ to describe the
assumed regularity conditions of $h$.

\textbf{Class ${R_{\beta}}$ :} $h(x)\in{R_{\beta}}$, if, with $\beta>0$ and
$x$ large enough, $h(x)$ can be represented as
\[
h(x)=x^{\beta}l(x),
\]
where $l(x)\in R_{0}$ and in $(\ref{3section102})$ $\epsilon(x)$ satisfies
\begin{equation}
\limsup_{x\rightarrow\infty}x|\epsilon^{\prime}(x)|<\infty,\qquad
\limsup_{x\rightarrow\infty}x^{2}|\epsilon^{^{\prime\prime}}(x)|<\infty.
\label{3section104}%
\end{equation}

\textbf{Class ${R_{\infty}}$ :} Further, $l\in\widetilde{R_{0}}$, if, in
$(\ref{3section102})$, $l(x)\rightarrow\infty$ as $x\rightarrow\infty$ and
\begin{equation}
\lim_{x\rightarrow\infty}\frac{x\epsilon^{\prime}(x)}{\epsilon(x)}%
=0,\qquad\lim_{x\rightarrow\infty}\frac{x^{2}\epsilon^{^{\prime\prime}}%
(x)}{\epsilon(x)}=0, \label{3section103}%
\end{equation}
and, for some $\eta\in(0,1/4)$
\begin{equation}
\liminf_{x\rightarrow\infty}x^{\eta}\epsilon(x)>0. \label{3section1030}%
\end{equation}
We say that $h\in{R_{\infty}}$ if $h$ is increasing and strictly monotone and
its inverse function $\psi$ defined through
\begin{equation}
\psi(u):=h^{\leftarrow}(u):=\inf\left\{  x:h(x)\geq u\right\}
\label{inverse de h}%
\end{equation}
belongs to $\widetilde{R_{0}}$.

Denote $\mathfrak{R:}={R_{\beta}}\cup{R_{\infty}}$. In fact, $\mathfrak{R} $
covers one large class of functions, although, ${R_{\beta}}$ and ${R_{\infty}%
}$ are only subsets of \emph{Regularly varying} and \emph{Rapidly varying}
functions, respectively.

\begin{rem}
The rôle of $(\ref{3section102})$ is to make $h(x)$ smooth enough. Under
$(\ref{3section102})$ the third order derivative of $h(x)$ exists, which is
necessary in order to use a Laplace methode for the asymptotic evaluation of
the moment generating function $\Phi(t)$ as $t\rightarrow\infty$, where
\[
\Phi(t)=\int_{0}^{\infty}e^{tx}p(x)dx=c\int_{0}^{\infty}\exp
\Big(K(x,t)+q(x)\Big)dx,\qquad t\in(0,\infty)
\]
in which
\[
K(x,t)=tx-g(x).
\]
If $h\in\mathfrak{R}$, $K(x,t)$ is concave with respect to $x$ and takes its
maximum at $\hat{x}=h^{\leftarrow}(t)$. \ The evaluation of $\Phi(t)$ for
large $t$ follows from an expansion of $K(x,t)$ in a neighborhood of $\hat
{x};$ this is Laplace's method. This expansion yields
\[
K(x,t)=K(\hat{x},t)-\frac{1}{2}h^{\prime}(\hat{x})\big(x-\hat{x}%
\big)^{2}-\frac{1}{6}h^{\prime\prime}(\hat{x})\big(x-\hat{x}\big)^{3}%
+\epsilon(x,t),
\]
where $\epsilon(x,t)$ is some error term. Conditions $(\ref{3section103})$
$(\ref{3section1030})$ and $(\ref{3section104})$ guarantee that $\epsilon
(x,t)$ goes to $0$ when $t$ tends to $\infty$ when $x$ belongs to some
neighborhood of $\hat{x}$.
\end{rem}

\begin{ex}
\textbf{Weibull Density. Let }$p$ be a Weibull density with shape parameter
$k>1$ and scale parameter $1$, namely%
\begin{align*}
p(x)  &  =kx^{k-1}\exp(-x^{k}),\qquad x\geq0\\
&  =k\exp\Big(-\big(x^{k}-(k-1)\log x\big)\Big).
\end{align*}
Take $g(x)=x^{k}-(k-1)\log x$ and $q(x)=0$. Then it holds
\[
h(x)=kx^{k-1}-\frac{k-1}{x}=x^{k-1}\big(k-\frac{k-1}{x^{k}}\big).
\]
Set $l(x)=k-(k-1)/x^{k},x\geq1$, then $(\ref{3section102})$ holds, namely,
\[
l(x)=\exp\Big(\int_{1}^{x}\frac{\epsilon(u)}{u}du\Big),\qquad x\geq1,
\]
with
\[
\epsilon(x)=\frac{k(k-1)}{kx^{k}-(k-1)}.
\]
The function $\epsilon$ is twice differentiable and goes to $0$ as
$x\rightarrow\infty$. Additionally, $\epsilon$ satisfies condition
$(\ref{3section104})$. Hence we have shown that $h\in R_{k-1}$.
\end{ex}

\begin{ex}
\textbf{A rapidly varying density.} Define $p$ through
\[
p(x)=c\exp(-e^{x-1}),\qquad x\geq0.
\]
Then $g(x)=h(x)=e^{x}$ and $q(x)=0$ for all non negative $x$. We show that
$h\in R_{\infty}$. It holds $\psi(x)=\log x+1$. Since $h(x)$ is increasing and
monotone, it remains to show that $\psi(x)\in\widetilde{R_{0}}$. When $x\geq
1$, $\psi(x)$ admits the representation of $(\ref{3section102})$ with
$\epsilon(x)=\log x+1$. Also conditions $(\ref{3section103})$ and
$(\ref{3section1030})$ are satisfied. Thus $h\in R_{\infty}$.
\end{ex}

Throughout the paper we use the following notation. When a r.v. $X$ has
density $p$ we write $p(X=x)$ instead of $p(x).$ This notation is useful when
changing measures. For example $\pi^{a}(X=x)$ is the density at point $x$ for
the variable $X$ generated under $\pi^{a}$, while $p(X=x)$ states for $X$
generated under $p.$ This avoids constant changes of notation.

\section{Conditional Density}

We inherit of the definition of the tilted density $\pi^{a}$ defined in
Section 1, and of the corresponding definitions of the functions $m$, $s^{2}$
and $\mu_{3}$. Because of (\ref{densityFunction}) and on the various
conditions on $g$ those functions are defined as $t\rightarrow\infty.$ The
following Theorem is basic for the proof of the remaining results.

\begin{thm}
\label{order of s} Let $p(x)$ be defined as in $(\ref{densityFunction})$ and
$h(x)\in\mathfrak{R}$. Denote by
\[
m(t)=\frac{d}{dt}\log\Phi(t),\quad\quad s^{2}(t)=\frac{d}{dt}m(t),\qquad
\mu_{3}(t)=\frac{d^{3}}{dt^{3}}\log\Phi(t),
\]
then with $\psi$ defined as in (\ref{inverse de h})it holds as $t\rightarrow
\infty$
\[
m(t)\sim\psi(t),\qquad s^{2}(t)\sim\psi^{\prime}(t),\qquad\mu_{3}(t)\sim
\frac{M_{6}-3}{2}\psi^{^{\prime\prime}}(t),
\]
where $M_{6}$ is the sixth order moment of standard normal distribution.
\end{thm}

The proof of this result relies on a series of Lemmas. Lemmas $(\ref{3lemma0}%
)$, $(\ref{3lemma01})$, $(\ref{3lemma02})$ and $(\ref{3lemma1})$ are used in
the proof. Lemma $(\ref{3lemma00})$ is instrumental for Lemma $(\ref{3lemma1}%
)$. The proof of Theorem $\ref{order of s}$ and these Lemmas are postponed to
Appendix. \newline

\begin{cor}
\label{3cor1} Let $p(x)$ be defined as in $(\ref{densityFunction})$ and
$h(x)\in\mathfrak{R}$. Then it holds as $t\rightarrow\infty$
\begin{align}
\label{3cor10g}\frac{\mu_{3}(t)}{s^{3}(t)}\longrightarrow0.
\end{align}

\end{cor}

Proof: Its proof relies on Theorem $2.1$ and is also put in Appendix.

\section{Edgeworth expansion under extreme normalizing factors}

With $\pi^{a_{n}}$ defined through
\[
\pi^{a_{n}}(x)=\frac{e^{tx}p(x)}{\Phi(t)},
\]
and $t$ determined by $a_{n}=m(t)$, define the normalized density of
$\pi^{a_{n}}$ by
\[
\bar{\pi}^{a_{n}}(x)=s_{n}\pi^{a_{n}}(s_{n}x+a_{n}),
\]
and denote the $n$-convolution of $\bar{\pi}^{a_{n}}(x)$ by $\bar{\pi}%
_{n}^{a_{n}}(x)$. Denote by $\rho_{n}$ the normalized density of
$n$-convolution $\bar{\pi}_{n}^{a_{n}}(x)$,
\[
\rho_{n}(x):=\sqrt{n}\bar{\pi}_{n}^{a_{n}}(\sqrt{n}x).
\]
The following result extends the local Edgeworth expansion of the distribution
of normalized sums of i.i.d. r;v's to the present context, where the summands
are generated under the density $\bar{\pi}^{a_{n}}$. Therefore the setting is
that of a triangular array of row wise independent summands; the fact that
$a_{n}\rightarrow\infty$ makes the situation unusual. We mainly adapt Feller's
proof (Chapiter 16, Theorem 2 $\cite{Feller}$).

\begin{thm}
\label{3theorem1} With the above notation, uniformly upon $x$ it holds
\[
\rho_{n}(x)=\phi(x)\Big(1+\frac{\mu_{3}}{6\sqrt{n}s^{3}}\big(x^{3}%
-3x\big)\Big)+o\Big(\frac{1}{\sqrt{n}}\Big).
\]
where $\phi(x)$ is standard normal density.
\end{thm}

Proof: \textbf{Step 1:} In this step, we will express the following formula
$G(x)$ by its Fourier transform. Let
\[
G(x):=\rho_{n}(x)-\phi(x)-\frac{\mu_{3}}{6\sqrt{n}s_{n}^{3}}\big(x^{3}%
-3x\big)\phi(x).
\]

From
\begin{align}
\label{3the010}\phi(x)=\frac{1}{2\pi}\int_{-\infty}^{\infty}e^{-i\tau
x}e^{-\frac{1}{2}\tau^{2}}d\tau,
\end{align}
it follows that
\begin{align}
\label{3the10p}\phi^{\prime\prime\prime}(x)=-\frac{1}{2\pi}\int_{-\infty
}^{\infty}(i\tau)^{3}e^{-i\tau x}e^{-\frac{1}{2}\tau^{2}}d\tau.
\end{align}
On the other hand
\[
\phi^{\prime\prime\prime}(x)=-(x^{3}-3x)\phi(x),
\]
which, together with $(\ref{3the10p})$, gives
\begin{equation}
(x^{3}-3x)\phi(x)=\frac{1}{2\pi}\int_{-\infty}^{\infty}(i\tau)^{3}e^{-i\tau
x}e^{-\frac{1}{2}\tau^{2}}d\tau. \label{3the11}%
\end{equation}

Let $\varphi^{a_{n}}(\tau)$ be the characteristic function (c.f) of $\bar{\pi
}^{a_{n}};$ the c.f of $\rho_{n}$ is $\big(\varphi^{a_{n}}(\tau/\sqrt
{n})\big)^{n}$. Hence it holds by Fourier inversion theorem
\begin{equation}
\rho_{n}(x)=\frac{1}{2\pi}\int_{-\infty}^{\infty}e^{-i\tau x}\big(\varphi
^{a_{n}}(\tau/\sqrt{n})\big)^{n}d\tau. \label{3the12}%
\end{equation}
Using $(\ref{3the010})$, $(\ref{3the11})$ and $(\ref{3the12})$, we have
\[
G(x)=\frac{1}{2\pi}\int_{-\infty}^{\infty}e^{-i\tau x}\Big(\big(\varphi
^{a_{n}}(\tau/\sqrt{n})\big)^{n}-e^{-\frac{1}{2}\tau^{2}}-\frac{\mu_{3}%
}{6\sqrt{n}s^{3}}(i\tau)^{3}e^{-\frac{1}{2}\tau^{2}}\Big)d\tau.
\]
Hence it holds
\begin{align}
&  \Big|\rho_{n}(x)-\phi(x)-\frac{\mu_{3}}{6\sqrt{n}s^{3}}\big(x^{3}%
-3x\big)\phi(x)\Big|\nonumber\label{3theo1000}\\
&  \leq\frac{1}{2\pi}\int_{-\infty}^{\infty}\Big|\big(\varphi^{a_{n}}%
(\tau/\sqrt{n})\big)^{n}-e^{-\frac{1}{2}\tau^{2}}-\frac{\mu_{3}}{6\sqrt
{n}s^{3}}(i\tau)^{3}e^{-\frac{1}{2}\tau^{2}}\Big|d\tau.
\end{align}

\textbf{Step 2:} In this step, we show that characteristic function
$\varphi^{a_{n}}$ of $\bar{\pi}^{a_{n}}(x)$ satisfies
\begin{align}
\label{3the130}\sup_{a_{n}\in\mathbb{R}^{+}}\int|\varphi^{a_{n}}(\tau
)|^{2}d\tau<\infty\qquad and\quad\sup_{a_{n}\in\mathbb{R}^{+},|\tau
|\geq\epsilon>0}|\varphi^{a_{n}}(\tau)|<1,
\end{align}
for any positive $\epsilon$ .

It is easy to verify that $r$-order ($r\geq1$) moment $\mu^{r}$ of $\pi
^{a_{n}}(x)$ satisfies
\[
\mu^{r}(t)=\frac{d^{r}\log\Phi(t)}{dt^{r}}\quad with\;t=m^{\leftarrow}%
(a_{n}),
\]
By Parseval identity
\begin{align}
\label{3the13}\int|\varphi^{a_{n}}(\tau)|^{2}d\tau=2\pi\int(\bar{\pi}^{a_{n}%
}(x))^{2}dx\leq2\pi\sup_{x\in\mathbb{R}}\bar{\pi}^{a_{n}}(x).
\end{align}
For the density function $p(x)$ in $(\ref{densityFunction})$, Theorem 5.4 of
Nagaev $\cite{Nagaev}$ states that the normalized conjugate density of $p(x)$,
namely, $\bar{\pi}^{a_{n}}(x)$ has the propriety
\[
\lim_{a_{n}\rightarrow\infty}\sup_{x\in\mathbb{R}}|\bar{\pi}^{a_{n}}%
(x)-\phi(x)|=0.
\]
Thus for arbitrary positive $\delta$, there exists some positive constant $M $
such that it holds
\[
\sup_{a_{n}\geq M}\sup_{x\in\mathbb{R}}|\bar{\pi}^{a_{n}}(x)-\phi
(x)|\leq\delta,
\]
which entails that $\sup_{a_{n}\geq M}\sup_{x\in\mathbb{R}}\bar{\pi}^{a_{n}%
}(x)<\infty$. When $a_{n}<M$, $\sup_{a_{n}<M}\sup_{x\in\mathbb{R}}\bar{\pi
}^{a_{n}}(x)<\infty;$ hence we have
\[
\sup_{a_{n}\in\mathbb{R}^{+}}\sup_{x\in\mathbb{R}}\bar{\pi}^{a_{n}}%
(x)<\infty,
\]
which, together with $(\ref{3the13})$, gives $(\ref{3the130})$. Furthermre,
$\varphi^{a_{n}}(\tau)$ is not periodic, hence the second inequality of
$(\ref{3the130})$ holds from Lemma $4$ (Chapiter $15$, section $1$) of
$\cite{Feller}$.

\textbf{Step 3:} In this step, we complete the proof by showing that when
$n\rightarrow\infty$
\begin{equation}
\int_{-\infty}^{\infty}\Big|\big(\varphi^{a_{n}}(\tau/\sqrt{n})\big)^{n}%
-e^{-\frac{1}{2}\tau^{2}}-\frac{\mu_{3}}{6\sqrt{n}s^{3}}(i\tau)^{3}%
e^{-\frac{1}{2}\tau^{2}}\Big|d\tau=o\Big(\frac{1}{\sqrt{n}}\Big).
\label{3theo100}%
\end{equation}

For arbitrarily positive sequence $a_{n}$ we have
\[
\sup_{a_{n}\in\mathbb{R}^{+}}\Big|\varphi^{a_{n}}(\tau)\Big|=\sup_{a_{n}%
\in\mathbb{R}^{+}}\Big|\int_{-\infty}^{\infty}e^{i\tau x}\bar{\pi}^{a_{n}%
}(x)dx\Big|\leq\sup_{a_{n}\in\mathbb{R}^{+}}\int_{-\infty}^{\infty
}\Big|e^{i\tau x}\bar{\pi}^{a_{n}}(x)\Big|dx=1.
\]
In addition, $\pi^{a_{n}}(x)$ is integrable, by Riemann-Lebesgue theorem, it
holds when $|\tau|\rightarrow\infty$
\[
\sup_{a_{n}\in\mathbb{R}^{+}}\Big|\varphi^{a_{n}}(\tau)\Big|\longrightarrow0.
\]
Thus for any strictly positive $\omega$, there exists some corresponding
$N_{\omega}$ such that if $|\tau|>\omega$, it holds
\begin{align}
\label{3the20}\sup_{a_{n}\in\mathbb{R}^{+}}\Big|\varphi^{a_{n}}(\tau
)\Big|<N_{\omega}<1.
\end{align}
We now turn to (\ref{3theo100}) which is splitted on $|\tau|>\omega\sqrt{n} $
and on $|\tau|\leq\omega\sqrt{n}$ .

It holds
\begin{align}
&  \sqrt{n}\int_{|\tau|>\omega\sqrt{n}}\Big|\big(\varphi^{a_{n}}(\tau/\sqrt
{n})\big)^{n}-e^{-\frac{1}{2}\tau^{2}}-\frac{\mu_{3}}{6\sqrt{n}s^{3}}%
(i\tau)^{3}e^{-\frac{1}{2}\tau^{2}}\Big|d\tau\nonumber\label{3the21}\\
&  \leq\sqrt{n}\int_{|\tau|>\omega\sqrt{n}}\Big|\big(\varphi^{a_{n}}%
(\tau/\sqrt{n})\big)\Big|^{n}d\tau+\sqrt{n}\int_{|\tau|>\omega\sqrt{n}%
}\Big|e^{-\frac{1}{2}\tau^{2}}+\frac{\mu_{3}}{6\sqrt{n}s^{3}}(i\tau
)^{3}e^{-\frac{1}{2}\tau^{2}}\Big|d\tau\nonumber\\
&  \leq\sqrt{n}N_{\omega}^{n-2}\int_{|\tau|>\omega\sqrt{n}}\Big|\big(\varphi
^{a_{n}}(\tau/\sqrt{n})\big)\Big|^{2}d\tau+\sqrt{n}\int_{|\tau|>\omega\sqrt
{n}}e^{-\frac{1}{2}\tau^{2}}\Big(1+\Big|\frac{\mu_{3}\tau^{3}}{6\sqrt{n}s^{3}%
}\Big|\Big)d\tau.
\end{align}
where the first term of the last line tends to $0$ when $n\rightarrow\infty$,
since
\begin{align}
&  \sqrt{n}N_{\omega}^{n-2}\int_{|\tau|>\omega\sqrt{n}}\Big|\big(\varphi
^{a_{n}}(\tau/\sqrt{n})\big)\Big|^{2}d\tau\nonumber\label{3the22}\\
&  =\exp\Big(\frac{1}{2}\log n+(n-2)\log N_{\omega}+\log\int_{|\tau
|>\omega\sqrt{n}}\Big|\big(\varphi^{a_{n}}(\tau/\sqrt{n})\big)\Big|^{2}%
d\tau\Big)\longrightarrow0,
\end{align}
where the last step holds from $(\ref{3the130})$ and $(\ref{3the20})$. As for
the second term of $(\ref{3the21})$, by Corollary $(\ref{3cor1})$, when
$n\rightarrow\infty$, we have $|\mu_{3}/s^{3}|\rightarrow0$. Hence it holds
when $n\rightarrow\infty$
\begin{align}
&  \sqrt{n}\int_{|\tau|>\omega\sqrt{n}}e^{-\frac{1}{2}\tau^{2}}%
\Big(1+\Big|\frac{\mu_{3}\tau^{3}}{6\sqrt{n}s^{3}}\Big|\Big)d\tau
\nonumber\label{3the23}\\
&  \leq\sqrt{n}\int_{|\tau|>\omega\sqrt{n}}e^{-\frac{1}{2}\tau^{2}}|\tau
|^{3}d\tau=\sqrt{n}\int_{|\tau|>\omega\sqrt{n}}\exp\Big\{-\frac{1}{2}\tau
^{2}+3\log|\tau|\Big\}d\tau\nonumber\\
&  =2\sqrt{n}\exp\big(-\omega^{2}n/2+o(\omega^{2}n/2)\big)\longrightarrow0,
\end{align}
where the second equality holds from, for example, Chapiter $4$ of
$\cite{Bingham}$. $(\ref{3the21})$, $(\ref{3the22})$ and $(\ref{3the23})$
implicate that, when $n\rightarrow\infty$
\begin{align}
\label{3theo10}\int_{|\tau|>\omega\sqrt{n}}\Big|\big(\varphi^{a_{n}}%
(\tau/\sqrt{n})\big)^{n}-e^{-\frac{1}{2}\tau^{2}}-\frac{\mu_{3}}{6\sqrt
{n}s^{3}}(i\tau)^{3}e^{-\frac{1}{2}\tau^{2}}\Big|d\tau=o\Big(\frac{1}{\sqrt
{n}}\Big).
\end{align}

If $|\tau|\leq\omega\sqrt{n}$, it holds
\begin{align}
&  \int_{|\tau|\leq\omega\sqrt{n}}\Big|\big(\varphi^{a_{n}}(\tau/\sqrt
{n})\big)^{n}-e^{-\frac{1}{2}\tau^{2}}-\frac{\mu_{3}}{6\sqrt{n}s^{3}}%
(i\tau)^{3}e^{-\frac{1}{2}\tau^{2}}\Big|d\tau\nonumber\label{3the24}\\
&  =\int_{|\tau|\leq\omega\sqrt{n}}e^{-\frac{1}{2}\tau^{2}}\Big|\big(\varphi
^{a_{n}}(\tau/\sqrt{n})\big)^{n}e^{\frac{1}{2}\tau^{2}}-1-\frac{\mu_{3}%
}{6\sqrt{n}s^{3}}(i\tau)^{3}\Big|d\tau\nonumber\\
&  =\int_{|\tau|\leq\omega\sqrt{n}}e^{-\frac{1}{2}\tau^{2}}\Big|\exp
\Big\{n\log\varphi^{a_{n}}(\tau/\sqrt{n})+{\frac{1}{2}\tau^{2}}\Big\}-1-\frac
{\mu_{3}}{6\sqrt{n}s^{3}}(i\tau)^{3}\Big|d\tau.
\end{align}
The integrand in the last display is bounded through
\begin{align}
\label{3the240}|e^{\alpha}-1-\beta|=|(e^{\alpha}-e^{\beta})+(e^{\beta}%
-1-\beta)|\leq(|\alpha-\beta|+\frac{1}{2}\beta^{2})e^{\gamma},
\end{align}
where $\gamma\geq\max(|\alpha|,|\beta|);$ this inequalityfollows replacing
$e^{\alpha},e^{\beta}$ by their power series, for real or complex
$\alpha,\beta$. Denote by
\[
\gamma(\tau)=\log\varphi^{a_{n}}(\tau)+{\frac{1}{2}\tau^{2}}.
\]
Since $\gamma^{\prime}(0)=\gamma^{\prime\prime}(0)=0$, the third order Taylor
expansion of $\gamma(\tau)$ at $\tau=0$ yields
\[
\gamma(\tau)=\gamma(0)+\gamma^{\prime}(0)\tau+\frac{1}{2}\gamma^{\prime\prime
}(0)\tau^{2}+\frac{1}{6}\gamma^{\prime\prime\prime}(\xi)\tau^{3}=\frac{1}%
{6}\gamma^{\prime\prime\prime}(\xi)\tau^{3},
\]
where $0<\xi<\tau$. Hence it holds
\[
\Big|\gamma(\tau)-\frac{\mu_{3}}{6s^{3}}(i\tau)^{3}\Big|=\Big|\gamma
^{\prime\prime\prime}(\xi)-\frac{\mu_{3}}{s_{n}^{3}}i^{3}\Big|\frac{|\tau
|^{3}}{6}.
\]
Here $\gamma^{\prime\prime\prime}$ is continuous; thus we can choose $\omega$
small enough such that $|\gamma^{\prime\prime\prime}(\xi)|<\rho$ for
$|\tau|<\omega$. Meanwhile, for $n$ large enough, according to Corollary
$(\ref{3cor1})$ , we have $|\mu_{3}/s^{3}|\rightarrow0$. Hence it holds for
$n$ large enough
\begin{align}
\label{3the25}\Big|\gamma(\tau)-\frac{\mu_{3}}{6s^{3}}(i\tau)^{3}%
\Big|\leq\Big(|\gamma^{\prime\prime\prime}(\xi)|+\rho\Big)\frac{|\tau|^{3}}%
{6}<\rho|\tau|^{3}.
\end{align}
Choose $\omega$ small enough, such that for $n$ large enough it holds for
$|\tau|<\omega$
\[
\Big|\frac{\mu_{3}}{6s^{3}}(i\tau)^{3}\Big|\leq\frac{1}{4}\tau^{2}%
,\qquad|\gamma(\tau)|\leq\frac{1}{4}\tau^{2}.
\]
For this choice of $\omega$, when $|\tau|<\omega$ we have
\begin{align}
\label{3the250}\max\Big(\Big|\frac{\mu_{3}}{6s^{3}}(i\tau)^{3}\Big|,|\gamma
(\tau)|\Big)\leq\frac{1}{4}\tau^{2}.
\end{align}
Replacing $\tau$ by $\tau/\sqrt{n}$, it holds for $|\tau|<\omega\sqrt{n}$
\begin{align}
&  \Big|n\log\varphi^{a_{n}}(\tau/\sqrt{n})+{\frac{1}{2}\tau^{2}}-\frac
{\mu_{3}}{6\sqrt{n}s^{3}}(i\tau)^{3}\Big|\nonumber\label{3the251}\\
&  =n\Big|\log\varphi^{a_{n}}(\tau/\sqrt{n})+{\frac{1}{2}\Big(\frac{\tau
}{\sqrt{n}}\Big)^{2}}-\frac{\mu_{3}}{6s^{3}}\Big(\frac{i\tau}{\sqrt{n}%
}\Big)^{3}\Big|\nonumber\\
&  =n\Big|\gamma\Big(\frac{\tau}{\sqrt{n}}\Big)-\frac{\mu_{3}}{6s^{3}%
}\Big(\frac{i\tau}{\sqrt{n}}\Big)^{3}\Big|<\frac{\rho|\tau|^{3}}{\sqrt{n}},
\end{align}
where the last inequality holds from $(\ref{3the25})$. In a similar way, with
$(\ref{3the250})$, it also holds for $|\tau|<\omega\sqrt{n}$
\begin{align}
&  \max\Big(\Big|n\log\varphi^{a_{n}}(\tau/\sqrt{n})+{\frac{1}{2}\tau^{2}%
}\Big|,\Big|\frac{\mu_{3}}{6\sqrt{n}s^{3}}(i\tau)^{3}%
\Big|\Big)\nonumber\label{3the252}\\
&  =n\max\Big(\Big|\gamma\Big(\frac{\tau}{\sqrt{n}}\Big)\Big|,\Big|\frac
{\mu_{3}}{6s^{3}}\Big(\frac{i\tau}{\sqrt{n}}\Big)^{3}\Big|\Big)\leq\frac{1}%
{4}\tau^{2}.
\end{align}

Apply $(\ref{3the240})$ to estimate the integrand of last line of
$(\ref{3the24})$, with the choice of $\omega$ in $(\ref{3the25})$ and
$(\ref{3the250})$, using $(\ref{3the251})$ and $(\ref{3the252})$ we have for
$|\tau|<\omega\sqrt{n}$
\begin{align*}
&  \Big|\exp\Big\{n\log\varphi^{a_{n}}(\tau/\sqrt{n})+{\frac{1}{2}\tau^{2}%
}\Big\}-1-\frac{\mu_{3}}{6\sqrt{n}s^{3}}(i\tau)^{3}\Big|\\
&  \leq\Big(\Big|n\log\varphi^{a_{n}}(\tau/\sqrt{n})+{\frac{1}{2}\tau^{2}%
}-\frac{\mu_{3}}{6\sqrt{n}s^{3}}(i\tau)^{3}\Big|+\frac{1}{2}\Big|\frac{\mu
_{3}}{6\sqrt{n}s^{3}}(i\tau)^{3}\Big|^{2}\Big)\\
&  \qquad\times\exp\Big[\max\Big(\Big|n\log\varphi^{a_{n}}(\tau/\sqrt
{n})+{\frac{1}{2}\tau^{2}}\Big|,\Big|\frac{\mu_{3}}{6\sqrt{n}s^{3}}(i\tau
)^{3}\Big|\Big)\Big]\\
&  \leq\Big(\frac{\rho|\tau|^{3}}{\sqrt{n}}+\frac{1}{2}\Big|\frac{\mu_{3}%
}{6\sqrt{n}s^{3}}(i\tau)^{3}\Big|^{2}\Big)\exp\Big(\frac{\tau^{2}}{4}\Big)\\
&  =\Big(\frac{\rho|\tau|^{3}}{\sqrt{n}}+\frac{\mu_{3}^{2}\tau^{6}}{72ns^{6}%
}\Big)\exp\Big(\frac{\tau^{2}}{4}\Big).
\end{align*}
Use this upper bound to $(\ref{3the24})$, we obtain
\begin{align*}
&  \int_{|\tau|\leq\omega\sqrt{n}}\Big|\big(\varphi^{a_{n}}(\tau/\sqrt
{n})\big)^{n}-e^{-\frac{1}{2}\tau^{2}}-\frac{\mu_{3}}{6\sqrt{n}s^{3}}%
(i\tau)^{3}e^{-\frac{1}{2}\tau^{2}}\Big|d\tau\\
&  \leq\int_{|\tau|\leq\omega\sqrt{n}}\exp\Big(-\frac{\tau^{2}}{4}%
\Big)\Big(\frac{\rho|\tau|^{3}}{\sqrt{n}}+\frac{\mu_{3}^{2}\tau^{6}}{72ns^{6}%
}\Big)d\tau\\
&  =\frac{\rho}{\sqrt{n}}\int_{|\tau|\leq\omega\sqrt{n}}\exp\Big(-\frac
{\tau^{2}}{4}\Big)|\tau|^{3}d\tau+\frac{\mu_{3}^{2}}{72ns^{6}}\int_{|\tau
|\leq\omega\sqrt{n}}\exp\Big(-\frac{\tau^{2}}{4}\Big)\tau^{6}d\tau,
\end{align*}
where both the first integral and the second integral are finite, and $\rho$
is arbitrarily small; additionally, by Corollary $(\ref{3cor1})$, ${\mu
_{3}^{2}}/{s^{6}}\rightarrow0$ when $n$ large enough, hence it holds when
$n\rightarrow\infty$
\begin{align}
\label{3theo11}\int_{|\tau|\leq\omega\sqrt{n}}\Big|\big(\varphi^{a_{n}}%
(\tau/\sqrt{n})\big)^{n}-e^{-\frac{1}{2}\tau^{2}}-\frac{\mu_{3}}{6\sqrt
{n}s^{3}}(i\tau)^{3}e^{-\frac{1}{2}\tau^{2}}\Big|d\tau=o\Big(\frac{1}{\sqrt
{n}}\Big).
\end{align}
Now $(\ref{3theo10})$ and $(\ref{3theo11})$ give $(\ref{3theo100})$. Further,
coming back to $(\ref{3theo1000})$, using $(\ref{3theo100})$, we obtain
\[
\Big|\bar{\pi}_{n}^{a_{n}}(x)-\phi(x)-\frac{\mu_{3}}{6\sqrt{n}s^{3}}%
\big(x^{3}-3x\big)\phi(x)\Big|=o\Big(\frac{1}{\sqrt{n}}\Big),
\]
which concludes the proof.

\section{Gibbs' conditional principles under extreme events}

We now explore Gibbs conditional principles under extreme events. The first
result is a pointwise approximation of the conditional density $p_{a_{n}%
}\left(  y_{1}^{k}\right)  $ on $\mathbb{R}^{k}$ for fixed $k.$ As a
by-product we also address the local approximation of $p_{A_{n}}$ where
\[
p_{A_{n}}\left(  y_{1}^{k}\right)  :=p\left(  \left.  X_{1}^{k}=y_{1}%
^{k}\right\vert S_{1}^{n}>na_{n}\right)  .
\]
However tese local approximations are of poor interest when comparing
$p_{a_{n}}$ to its approximation.\ 

We consider the case $k=1.$ For $Y_{1}^{n}$ a random vector with density
$p_{a_{n}}$ we first provide a density $g_{a_{n}}$ on $\mathbb{R}$ such that
\[
p_{a_{n}}\left(  Y_{1}\right)  =g_{a_{n}}\left(  Y_{1}\right)  \left(
1+R_{n}\right)
\]
where $R_{n}$ is a function of the vector $Y_{1}^{n}$ which goes to $0$ as $n$
tends to infinity. The above statement may also be written as
\begin{equation}
p_{a_{n}}\left(  y_{1}\right)  =g_{a_{n}}\left(  y_{1}\right)  \left(
1+o_{P_{a_{n}}}(1)\right)  \label{approxSurRandom}%
\end{equation}
where $P_{a_{n}}$ is the joint probability measure of the vector $Y_{1}^{n}$
under the condition $\left(  S_{1}^{n}=na_{n}\right)  .$ This statement is of
a different nature with respect to the above one, since it amounts to prove
the approximation on typical realisations under the conditional sampling
scheme. We will deduce from (\ref{approxSurRandom}) that the $L^{1}$ distance
between $p_{a_{n}}$ and $g_{a_{n}}$ goes to $0$ as $n$ tends to infinity. It
would be interesting to extend these results to the case when $k=k_{n}$ is
close to $n$, as done in \cite{Bronia} in all cases from the CLT to the LDP
ranges. The extreme deviation case is more envolved, which led us to restrict
this study to the case when $k=1$ (or $k$ fixed, similarly).

\subsection{\bigskip A local result in $\mathbb{R}^{k}$}

Fix $y_{1}^{k}:=\left(  y_{1},..,y_{k}\right)  $ in $\mathbb{R}^{k}$ and
define $s_{i}^{j}:=y_{i}+..+y_{j}$ for $1\leq i<j\leq k.$

Define $t_{i}$ through
\begin{equation}
m(t_{i}):=\frac{na_{n}-s_{1}^{i}}{n-i}. \label{3lfd01}%
\end{equation}
For the sake of brevity, we write $m_{i}$ instead of $m(t_{i})$, and define
$s_{i}^{2}:=s^{2}(t_{i})$. We have the following conditional density.

Consider the following condition%
\begin{equation}
\lim_{n\rightarrow\infty}\frac{\psi(t_{k})^{2}}{\sqrt{n\psi^{\prime}(t_{k})}%
}=0, \label{croissance de a}%
\end{equation}
which can be seen as a growth condition on $a_{n}$, avoiding too large
increases of this sequence.

For $0\le i \le k-1<n $, define $z_{i}$ through
\begin{align*}
z_{i}=\frac{m_{i}-y_{i+1}}{s_{i} \sqrt{n-i-1}}.
\end{align*}

\begin{lem}
\label{3lemma z} Assume that $p(x)$ satisfies $(\ref{densityFunction})$ and
$h(x)\in\mathfrak{R}$. Let $t_{i}$ is defined in $(\ref{3lfd01})$. Assume that
$a_{n}\rightarrow\infty$ as $n\rightarrow\infty$ and that
(\ref{croissance de a}) holds. then it holds as $a_{n}\rightarrow\infty$
\[
\lim_{n\rightarrow\infty}\sup_{0\leq i\leq k-1}z_{i}=0.\qquad
\]

\end{lem}

Proof: When $n\rightarrow\infty$, it holds
\[
z_{i}\sim{m_{i}}/{s_{i}\sqrt{n-i-1}}\sim m_{i}/(s_{i}\sqrt{n}).
\]
From Theorem $\ref{order of s}$, it holds $m(t)\sim{\psi(t)}$ and
$s(t)\sim\sqrt{\psi^{^{\prime}}(t)}$. Hence we have
\begin{align}
\label{3felff}z_{i}\sim\frac{\psi(t_{i})}{\sqrt{n\psi^{^{\prime}}(t_{i})}}.
\end{align}
By $(\ref{3lfd01})$, $m_{i}\sim m_{k}$ as $n\rightarrow\infty$. Consider
$m_{k} \sim\psi(t_{k})$. Then it holds
\begin{align*}
m_{i} \sim\psi(t_{k}).
\end{align*}
In addition, $m_{i} \sim\psi(t_{i})$ by Theorem $\ref{order of s}$, this
implies it holds
\begin{align}
\label{3fel}\psi(t_{i}) \sim\psi(t_{k}).
\end{align}

\textbf{Case 1:} if $h(x)\in R_{\beta}$. We have $h(x)=x^{\beta}l_{0}(x) ,
l_{0}(x)\in R_{0},\beta>0$. Hence
\begin{align*}
h^{^{\prime}}(x)=x^{\beta-1}l_{0}(x)\big(\beta+\epsilon(x)\big),
\end{align*}
set $x=\psi(t)$, we get
\begin{align}
\label{3fel0}h^{^{\prime}}\big(\psi(t)\big)=\big(\psi(t)\big)^{\beta-1}%
l_{0}\big(\psi(t)\big)\big(\beta+\epsilon\big(\psi(t)\big)\big).
\end{align}
Notice that it holds $\psi^{^{\prime}}(t)=1/h^{^{\prime}}\big(\psi(t)\big)$,
combine $(\ref{3fel})$ with $(\ref{3fel0})$, we obtain
\begin{align*}
\frac{\psi^{\prime}(t_{i})}{\psi^{\prime}(t_{k})}=\frac{h^{^{\prime}}%
\big(\psi(t_{k})\big)}{h^{^{\prime}}\big(\psi(t_{i})\big)} =\frac
{\big(\psi(t_{k})\big)^{\beta-1}l_{0}\big(\psi(t_{k})\big)\big(\beta
+\epsilon\big(\psi(t_{k})\big)\big)} {\big(\psi(t_{i})\big)^{\beta-1}%
l_{0}\big(\psi(t_{i})\big)\big(\beta+\epsilon\big(\psi(t_{i})\big)\big)}%
\longrightarrow1,
\end{align*}
where we use the slowly varying propriety of $l_{0}$. Thus it holds
\begin{align*}
\psi^{\prime}(t_{i}) \sim\psi^{\prime}(t_{k}),
\end{align*}
which, together with $(\ref{3fel})$, is put into $(\ref{3felff})$ to yield
\begin{align*}
z_{i} \sim\frac{\psi(t_{k})}{\sqrt{n\psi^{^{\prime}}(t_{k})}}.
\end{align*}
Hence we have under condition $(\ref{croissance de a})$
\begin{align*}
z_{i}^{2} \sim\frac{\psi(t_{k})^{2}}{{n\psi^{^{\prime}}(t_{k})}}=\frac
{\psi(t_{k})^{2}}{{\sqrt{n}\psi^{^{\prime}}(t_{k})}}\frac{1}{\sqrt{n}%
}=o\Big(\frac{1}{\sqrt{n}}\Big),
\end{align*}
which implies further $z_{i}\rightarrow0$. Note that the final step is used in
order to relax the strength of the growth condition on $a_{n}.$

\textbf{Case 2:} if $h(x)\in R_{\infty}$. By $(\ref{3lfd01})$, it holds
$m(t_{k})\geq m(t_{i})$ as $n\rightarrow\infty$. Since the function
$t\rightarrow m(t)$ is increasing, we have
\[
t_{i}\leq t_{k}.
\]
The function $t\rightarrow\psi^{^{\prime}}(t)$ is decreasing, since
\[
\psi^{^{\prime\prime}}(t)=-\frac{\psi(t)}{t^{2}}\epsilon
(t)\big(1+o(1)\big)<0\qquad as\quad t\rightarrow\infty.
\]
Therefore it holds as $n\rightarrow\infty$
\[
\psi^{\prime}(t_{i})\geq\psi^{\prime}(t_{k}),
\]
which, combined with $(\ref{3felff})$ and $(\ref{3fel})$, yields
\[
z_{i}\sim\frac{\psi(t_{k})}{\sqrt{n\psi^{^{\prime}}(t_{i})}}\leq\frac
{2\psi(t_{k})}{\sqrt{n\psi^{^{\prime}}(t_{k})}},
\]
hence we have
\[
z_{i}^{2}\leq\frac{4\psi(t_{k})^{2}}{{n\psi^{^{\prime}}(t_{k})}}=\frac
{4\psi(t_{k})^{2}}{{\sqrt{n}\psi^{^{\prime}}(t_{k})}}\frac{1}{\sqrt{n}%
}=o\Big(\frac{1}{\sqrt{n}}\Big),
\]
where the last step holds from condition $(\ref{croissance de a})$. Further it
holds $z_{i}\rightarrow0$.

\begin{thm}
\label{point conditional density} With the above notation and hypotheses,
assuming (\ref{croissance de a}), it holds
\[
p_{a_{n}}(y_{1}^{k})=p(X_{1}^{k}=y_{1}^{k}|S_{1}^{n}=na_{n})=g_{m}(y_{1}%
^{k})\Big(1+o(1)\Big).
\]
with
\[
g_{m}(y_{1}^{k})=\prod_{i=0}^{k-1}\Big(\pi^{m_{i}}(X_{i+1}=y_{i+1})\Big).
\]

\end{thm}

Proof:

Using Bayes formula,
\begin{align}
p_{a_{n}}\left(  y_{1}^{k}\right)  :=p(X_{1}^{k}=y_{1}^{k}|S_{1}^{n}=na_{n})
&  =p(X_{1}=y_{1}|S_{1}^{n}=na_{n})\prod_{i=1}^{k-1}p(X_{i+1}=y_{i+1}%
|X_{1}^{i}=y_{1}^{i},S_{1}^{n}=na_{n})\nonumber\label{bayes formula01}\\
&  =\prod_{i=0}^{k-1}p(X_{i+1}=y_{i+1}|S_{i+1}^{n}=na_{n}-s_{1}^{i}).
\end{align}
We make use of the following invariance property:For all $y_{1}^{k}$ and all
$\alpha>0$%
\[
p(X_{i+1}=y_{i+1}|X_{1}^{i}=y_{1}^{i},S_{1}^{n}=na_{n})=\pi^{\alpha}%
(X_{i+1}=y_{i+1}|X_{1}^{i}=y_{1}^{i},S_{1}^{n}=na_{n})
\]
where on the LHS, the r.v's $X_{1}^{i}$ are sampled i.i.d. under $p$ and on
the RHS, sampled i.i.d. under $\pi^{\alpha}.Itthusholds$%

\begin{align}
&  p(X_{i+1}=y_{i+1}|S_{i+1}^{n}=na_{n}-S_{1}^{i})=\pi^{m_{i}}(X_{i+1}%
=y_{i+1}|S_{i+1}^{n}=na_{n}-s_{1}^{i})\nonumber\label{bayes formula}\\
&  =\pi^{m_{i}}(X_{i+1}=y_{i+1})\frac{\pi^{m_{i}}(S_{i+2}^{n}=na_{n}%
-s_{1}^{i+1})}{\pi^{m_{i}}(S_{i+1}^{n}=na_{n}-s_{1}^{i})}\nonumber\\
&  =\frac{\sqrt{n-i}}{\sqrt{n-i-1}}\pi^{m_{i}}(X_{i+1}=y_{i+1})\frac
{\widetilde{\pi_{n-i-1}}(\frac{m_{i}-y_{i+1}}{s_{i}\sqrt{n-i-1}})}%
{\widetilde{\pi_{n-i}}(0)},
\end{align}
where $\widetilde{\pi_{n-i-1}}$ is the normalized density of $S_{i+2}^{n}$
under i.i.d. sampling under $\pi^{m_{i}};$correspondingly, $\widetilde
{\pi_{n-i}}$ is the normalized density of $S_{i+1}^{n}$ under the same
sampling. Note that a r.v. with density $\pi^{mi}$ has expectation $m_{i}$ and
variance $s_{i}^{2}$.

Write $z_{i}=\frac{m_{i}-y_{i+1}}{s_{i}\sqrt{n-i-1}}$, and perform a
third-order Edgeworth expansion of $\widetilde{\pi_{n-i-1}}(z_{i})$, using
Theorem $\ref{3theorem1}$. It follows
\begin{align}
\label{edgeworth exte}\widetilde{\pi_{n-i-1}}(z_{i})=\phi(z_{i})\Big(1+\frac
{\mu_{3}^{i}}{6s_{i}^{3}\sqrt{n-1}}(z_{i}^{3}-3z_{i})\Big)+o\Big(\frac
{1}{\sqrt{n}}\Big),
\end{align}
The approximation of $\widetilde{\pi_{n-i}}(0)$ is obtained from
$(\ref{edgeworth exte})$
\begin{align}
\label{pi Z=0}\widetilde{\pi_{n-i}}(0)=\phi(0)\Big(1+o\big(\frac{1}{\sqrt{n}%
}\big)\Big).
\end{align}
Put $(\ref{edgeworth exte})$ and $(\ref{pi Z=0})$ into $(\ref{bayes formula})
$ to obtain
\begin{align}
&  p(X_{i+1}=y_{i+1}|S_{i+1}^{n}=na_{n}-S_{1}^{i}%
)\nonumber\label{3conditionaldensity01}\\
&  =\frac{\sqrt{n-i}}{\sqrt{n-i-1}}\pi^{m_{i}}(X_{i+1}=y_{i+1})\frac
{\phi(z_{i})}{\phi(0)}\Big[1+\frac{\mu_{3}^{i}}{6s_{i}^{3}\sqrt{n-1}}%
(z_{i}^{3}-3z_{i})+o\Big(\frac{1}{\sqrt{n}}\Big)\Big]\nonumber\\
&  =\frac{\sqrt{2\pi(n-i)}}{\sqrt{n-i-1}}\pi^{m_{i}}(X_{i+1}=y_{i+1}%
){\phi(z_{i})}\big(1+R_{n}+o(1/\sqrt{n})\big),
\end{align}
where
\[
R_{n}=\frac{\mu_{3}^{i}}{6s_{i}^{3}\sqrt{n-1}}(z_{i}^{3}-3z_{i}).
\]

Under condition $(\ref{croissance de a})$, using Lemma $\ref{3lemma z}$, it
holds $z_{i}\rightarrow0$ as $a_{n}\rightarrow\infty$, and under Corollary
$(\ref{3cor1})$, $\mu_{3}^{i}/s_{i}^{3}\rightarrow0.$ This yields
\[
R_{n}=o\big(1/\sqrt{n}\big),
\]
which, combined with $(\ref{3conditionaldensity01})$, gives
\begin{align*}
&  p(X_{i+1}=y_{i+1}|s_{i+1}^{n}=na_{n}-S_{1}^{i})=\frac{\sqrt{2\pi(n-i)}%
}{\sqrt{n-i-1}}\pi^{m_{i}}(X_{i+1}=y_{i+1}){\phi(z_{i})}\big(1+o(1/\sqrt
{n})\big)\\
&  =\frac{\sqrt{n-i}}{\sqrt{n-i-1}}\pi^{m_{i}}(X_{i+1}=y_{i+1}){\big(1-z_{i}%
^{2}/2+o(z_{i}^{2})\big)}\big(1+o(1/\sqrt{n})\big),
\end{align*}
where we use one Taylor expansion in second equality. Using once more Lemma
$\ref{3lemma z}$, under conditions $(\ref{croissance de a})$, we have as
$a_{n}\rightarrow\infty$
\[
z_{i}^{2}=o(1/\sqrt{n}),
\]
hence we get
\[
p(X_{i+1}=y_{i+1}|S_{i+1}^{n}=na_{n}-s_{1}^{i})=\frac{\sqrt{n-i}}{\sqrt
{n-i-1}}\pi^{m_{i}}(X_{i+1}=y_{i+1})\big(1+o(1/\sqrt{n})\big),
\]
which together with $(\ref{bayes formula01})$ yields
\begin{align*}
p(X_{1}^{k}=y_{1}^{k}|S_{1}^{n}=na_{n})  &  =\prod_{i=0}^{k-1}\Big(\frac
{\sqrt{n-i}}{\sqrt{n-i-1}}\pi^{m_{i}}(X_{i+1}=y_{i+1})\big(1+o(1/\sqrt
{n})\big)\Big)\\
& =\prod_{i=0}^{k-1}\Big(\pi^{m_{i}}(X_{i+1}=y_{i+1}) \Big)\prod_{i=0}%
^{k-1}\Big(\frac{\sqrt{ n-i}}{\sqrt{n-i-1}}\Big)\prod_{i=0}^{k-1}%
\Big( 1+o\big(\frac{1}{\sqrt{n}}\big)\Big)\\
& =\Big( 1+o\big(\frac{1}{\sqrt{n}}\big)\Big)\prod_{i=0}^{k-1}\Big(\pi^{m_{i}%
}(X_{i+1}=y_{i+1}) \Big),
\end{align*}
The proof is completed.

\bigskip Define $t$ through $m(t)=a_{n}$, replace condition
$(\ref{croissance de a})$ by
\begin{equation}
\lim_{n\rightarrow\infty}\frac{\psi(t)^{2}}{\sqrt{n\psi^{\prime}(t)}}=0,
\label{croissance de a n}%
\end{equation}
then for fixed $k$, an equivalent statement is

\begin{thm}
\label{point conditional density e} Under the same hypotheses as in the
previous Theorem
\[
p_{a_{n}}(y_{1}^{k})=p(X_{1}^{k}=y_{1}^{k}|S_{1}^{n}=na_{n})=g_{a_{n}}%
(y_{1}^{k})\Big(1+o\big(\frac{1}{\sqrt{n}}\big)\Big).
\]
with
\[
g_{a_{n}}(y_{1}^{k})=\prod_{i=1}^{k}\Big(\pi^{a_{n}}(X_{i}=y_{i})\Big).
\]

\end{thm}

Proof:

Using the notations of Theorem $\ref{point conditional density}$, by
$(\ref{bayes formula01})$, we obtain%

\begin{align}
\label{bayes formula01 e}p(X_{1}^{k}=y_{1}^{k}|S_{1}^{n}=na_{n})=\prod
_{i=0}^{k-1}p(X_{i+1}=y_{i+1}|S_{i+1}^{n}=na_{n}-S_{1}^{i}).
\end{align}
$(\ref{bayes formula})$ is replaced by
\begin{equation}
p(X_{i+1}=y_{i+1}|S_{i+1}^{n}=na_{n}-S_{1}^{i})=\frac{\sqrt{n-i}}{\sqrt
{n-i-1}}\pi^{a_{n}}(X_{i+1}=y_{i+1})\frac{\widetilde{\pi_{n-i-1}^{a_{n}}%
}(\frac{(i+1)a_{n}-S_{1}^{i+1}}{s\sqrt{n-i-1}})}{\widetilde{\pi_{n-i}^{a_{n}}%
}\big(\frac{ia_{n}-S_{1}^{i}}{s\sqrt{n-i}}\big)}, \label{bayes formula e}%
\end{equation}
where $\widetilde{\pi_{n-i-1}^{a_{n}}}({(i+1)a_{n}-y_{i+1}}/{s_{i}\sqrt
{n-i-1}})$ is the normalized density of $\pi^{a_{n}}(S_{i+2}^{n}=na_{n}%
-S_{1}^{i+1})$, and $\pi^{a_{n}}$ has the expectation $a_{n}$ and variance
$s$. Correspondingly, $\widetilde{\pi_{n-i}^{a_{n}}}\big(({ia_{n}-S_{1}^{i}%
})/{s\sqrt{n-i}}\big)$ is the normalized density of $\pi^{a_{n}}(S_{i+1}%
^{n}=na_{n}-S_{1}^{i})$.

Write $z_{i}=\frac{(i+1)a_{n}-S_{1}^{i+1}}{s\sqrt{n-i-1}}$, by Theorem
$\ref{3theorem1}$ one three-order Edgeworth expansion yields
\begin{align}
\label{edgeworth exte e}\widetilde{\pi_{n-i-1}^{a_{n}}}(z_{i})=\phi
(z_{i})\big(1+R_{n}^{i}\big)+o\Big(\frac{1}{\sqrt{n}}\Big),
\end{align}
where
\[
R_{n}^{i}=\frac{\mu_{3}}{6s^{3}\sqrt{n-1}}(z_{i}^{3}-3z_{i}).
\]
Set $i=i-1$, the approximation of $\widetilde{\pi_{n-i}^{a_{n}}}$ is obtained
from $(\ref{edgeworth exte e})$
\begin{align}
\label{pi Z=0 e}\widetilde{\pi_{n-i}}(z_{i-1})=\phi(z_{i-1})\Big(1+R_{n}%
^{i+1}\Big)+o\Big(\frac{1}{\sqrt{n}}\Big).
\end{align}
When $a_{n}\rightarrow\infty$, using Theorem $\ref{order of s}$, it holds
\begin{align}
\sup_{0\leq i\leq k-1}z_{i}^{2}\sim\frac{(i+1)^{2}a_{n}^{2}}{s^{2}{n}}  &
\leq\frac{2k^{2}a_{n}^{2}}{s^{2}{n}}=\frac{2k^{2}(m(t))^{2}}{s^{2}{n}%
}\nonumber\label{3fqmq}\\
&  \sim\frac{2k^{2}(\psi(t))^{2}}{\psi^{\prime}(t){n}}=\frac{2k^{2}%
(\psi(t))^{2}}{\sqrt{n}\psi^{\prime}(t)}\frac{1}{\sqrt{n}}=o\Big(\frac
{1}{\sqrt{n}}\Big),
\end{align}
where last step holds under condition $(\ref{croissance de a n})$. Hence it
holds $z_{i}\rightarrow0$ uniformly in $i$ as $a_{n}\rightarrow\infty$, and by
Corollary $(\ref{3cor1})$, $\mu_{3}/s^{3}\rightarrow0$, then it follows
\[
R_{n}^{i}=o\big(1/\sqrt{n}\big)\qquad R_{n}^{i+1}=o\big(1/\sqrt{n}\big),
\]
then put $(\ref{edgeworth exte e})$ and $(\ref{pi Z=0 e})$ into
$(\ref{bayes formula e})$, we obtain
\begin{align*}
&  p(X_{i+1}=y_{i+1}|S_{i+1}^{n}=na_{n}-S_{1}^{i})=\frac{\sqrt{n-i}}%
{\sqrt{n-i-1}}\pi^{a_{n}}(X_{i+1}=y_{i+1})\frac{\phi(z_{i})}{\phi(z_{i-1}%
)}\big(1+o(1/\sqrt{n})\big)\\
&  =\frac{\sqrt{n-i}}{\sqrt{n-i-1}}\pi^{m_{i}}(X_{i+1}=y_{i+1}){\big(1-(z_{i}%
^{2}-z_{i-1}^{2})/2+o(z_{i}^{2}-z_{i-1}^{2})\big)}\big(1+o(1/\sqrt{n})\big),
\end{align*}
where we use one Taylor expansion in second equality. Using $(\ref{3fqmq})$,
we have as $a_{n}\rightarrow\infty$
\[
|z_{i}^{2}-z_{i-1}^{2}|=o(1/\sqrt{n}),
\]
hence we get
\[
p(X_{i+1}=y_{i+1}|S_{i+1}^{n}=na_{n}-S_{1}^{i})=\frac{\sqrt{n-i}}{\sqrt
{n-i-1}}\pi^{a_{n}}(X_{i+1}=y_{i+1})\big(1+o(1/\sqrt{n})\big),
\]
which together with $(\ref{bayes formula01 e})$ yields
\begin{align}
p(X_{1}^{k}=y_{1}^{k}|S_{1}^{n}=na_{n})  &  =\prod_{i=0}^{k-1}\Big(\pi^{a_{n}%
}(X_{i+1}=y_{i+1})\sqrt{\frac{n}{n-k}}\Big)\prod_{i=0}^{k-1}\Big(1+o\big(\frac
{1}{\sqrt{n}}\big)\Big)\nonumber\label{3pfl00 e}\\
&  =\Big(1+o\big(\frac{1}{\sqrt{n}}\big)\Big)\prod_{i=0}^{k-1}\Big(\pi^{a_{n}%
}(X_{i+1}=y_{i+1})\Big).
\end{align}
This completes the proof. \bigskip\bigskip

\begin{rem}
The above result shows that asymptotically the point condition $\left(
S_{1}^{n}=na_{n}\right)  $ leaves blocks of $k$ of the $X_{i}^{\prime}s$
independent. Obviously this property does not hold for large values of $k,$
close to $n.$ A similar statement holds in the LDP range, conditioning either
on $\left(  S_{1}^{n}=na\right)  $ (see Diaconis and Friedman 1988)), or on
$\left(  S_{1}^{n}\geq na\right)  $; see Csiszar 1984 for a general statement
on asymptotic conditional independence.
\end{rem}

Using the same proof as in Theorem $(\ref{point conditional density e})$, we
obtain the following corollary.

\begin{cor}
\label{3cor2} It holds
\[
p_{a}(y_{1}^{k})=p(X_{1}^{k}=y_{1}^{k}|S_{1}^{n}=na_{n})=g_{a}(y_{1}%
^{k})\Big(1+o\big(\frac{k}{\sqrt{n}}\big)\Big).
\]
with
\[
g_{a}(y_{1}^{k})=\prod_{i=1}^{k}\Big(\pi^{a}(X_{i}=y_{i})\Big).
\]

\end{cor}

\bigskip

\subsection{Strenghtening of the local Gibbs conditional principle}

We now turn to a stronger approximation of $p_{a_{n}}.$ Consider $Y_{1}^{n}$
with density $p_{a_{n}}$ and the resulting random variable p$_{a_{n}}\left(
Y_{1}\right)  .$ We prove the following result

\begin{thm}
\label{ThmLocalProba}With all the above notation and hypotheses it holds%
\[
p_{a_{n}}\left(  Y_{1}\right)  =g_{a_{n}}\left(  Y_{1}\right)  \left(
1+R_{n}\right)
\]
where
\[
g_{a_{n}}=\pi^{a_{n}}%
\]
the tilted density at point $a_{n}$ , and where $R_{n}$ is a function of
$Y_{1}^{n}$ such that $P_{a_{n}}\left(  \left\vert R_{n}\right\vert
>\delta\sqrt{n}\right)  \rightarrow0$ as $n\rightarrow\infty$ for any positive
$\delta.$
\end{thm}

\bigskip

This result is of much greater relevance than the previous ones. Indeed under
$P_{a_{n}}$ the r.v. $Y_{1}$ may take large values. At the contrary simple
approximation of $p_{a_{n}}$ by $g_{a_{n}}$ on $\mathbb{R}_{+}$ only provides
some knowledge on $p_{a_{n}}$ on sets with smaller and smaller probability
under $p_{a_{n}}.$ Also it will be proved that as a consequence of the above
result, the $L^{1}$ norm between $p_{a_{n}}$ and $g_{a_{n}}$ goes to $0$ as
$n\rightarrow\infty$, a result out of reach through the aforementioned results.

In order to adapt the proof of Theorem *** to the present setting it is
necessary to get some insight on the plausible values of $Y_{1}$ under
$P_{a_{n}}.$ It holds

\begin{lem}
Under $P_{a_{n}}$ it holds%
\[
Y_{1}=O_{P_{a_{n}}}\left(  a_{n}\right)
\]

\end{lem}

Proof: This is a consequence of Markov Inequality:%

\[
P\left(  \left.  Y_{1}>u\right\vert S_{1}^{n}=na_{n}\right)  \leq
\frac{E\left(  \left.  Y_{1}\right\vert S_{1}^{n}=na_{n}\right)  }{u}%
=\frac{a_{n}}{u}%
\]
\bigskip which goes to $0$ for all $u=u_{n}$ such that lim$_{n\rightarrow
\infty}u_{n}/a_{n}=\infty.$

\bigskip

We now turn back to the proof of our result, replacing $y_{1}^{k}$ by $Y_{1}$
in (\ref{bayes formula e}).

It holds%
\[
P\left(  \left.  X_{1}=Y_{1}\right\vert S_{1}^{n}=na_{n}\right)
=P(X_{1}=Y_{1})\frac{P\left(  S_{2}^{n}=na_{n}-Y_{1}\right)  }{P\left(
S_{1}^{n}=na_{n}\right)  }%
\]
in which the tilting substitution of measures is performed, with tilting
density $\pi^{a_{n}}$, followed by normalization. Now if the growth condition
(\ref{croissance de a}) holds, namely%
\[
\lim_{n\rightarrow\infty}\frac{\psi(t)}{\sqrt{n\psi^{\prime}(t)}}=0
\]
with $m(t)=a_{n}$ $\ $\ it follows that
\[
P\left(  \left.  X_{1}=Y_{1}\right\vert S_{1}^{n}=na_{n}\right)  =\pi^{a_{n}%
}\left(  Y_{1}\right)  \left(  1+R_{n}\right)
\]
as claimed where the order of magnitude of $R_{n}$ is $o_{P_{a_{n}}}\left(
1/\sqrt{n}\right)  $. We have proved Theorem \ref{ThmLocalProba}.

Denote the conditional probabilities by $P_{a_{n}}$ and $G_{a_{n}}$ which
correspond to the density functions $p_{a_{n}}$ and $g_{a_{n}}$, respectively.

\subsection{Gibbs principle in variation norm}

We now consider the approximation of $P_{a_{n}}$ by $G_{a_{n}}$ in variation norm.

The main ingredient is the fact that in the present setting approximation of
$p_{a_{n}}$ by $g_{a_{n}}$ in probability plus some rate implies approximation
of the corresponding measures in variation norm. This approach has been
developped in Broniatowski and Caron (2012); we state a first lemma which
states that wether two densities are equivalent in probability with small
relative error when measured according to the first one, then the same holds
under the sampling of the second.

Let $\mathfrak{R}_{n}$ and $\mathfrak{S}_{n}$ denote two p.m's on
$\mathbb{R}^{n}$ with respective densities $\mathfrak{r}_{n}$ and
$\mathfrak{s}_{n}.$

\begin{lem}
\label{Lemma:commute_from_p_n_to_g_n} Suppose that for some sequence
$\varepsilon_{n}$ which tends to $0$ as $n$ tends to infinity%
\begin{equation}
\mathfrak{r}_{n}\left(  Y_{1}^{n}\right)  =\mathfrak{s}_{n}\left(  Y_{1}%
^{n}\right)  \left(  1+o_{\mathfrak{R}_{n}}(\varepsilon_{n})\right)
\label{p_n equiv g_n under p_n}%
\end{equation}
as $n$ tends to $\infty.$ Then
\begin{equation}
\mathfrak{s}_{n}\left(  Y_{1}^{n}\right)  =\mathfrak{r}_{n}\left(  Y_{1}%
^{n}\right)  \left(  1+o_{\mathfrak{S}_{n}}(\varepsilon_{n})\right)
.\label{g_n equiv p_n under g_n}%
\end{equation}

\end{lem}

\begin{proof}
Denote
\begin{equation*}
A_{n,\varepsilon_{n}}:=\left\{ y_{1}^{n}:(1-\varepsilon_{n})\mathfrak{s}%
_{n}\left( y_{1}^{n}\right) \leq\mathfrak{r}_{n}\left( y_{1}^{n}\right) \leq%
\mathfrak{s}_{n}\left( y_{1}^{n}\right) (1+\varepsilon_{n})\right\} .
\end{equation*}
It holds for all positive $\delta$%
\begin{equation*}
\lim_{n\rightarrow\infty}\mathfrak{R}_{n}\left(
A_{n,\delta\varepsilon_{n}}\right) =1.
\end{equation*}
Write
\begin{equation*}
\mathfrak{R}_{n}\left( A_{n,\delta\varepsilon_{n}}\right) =\int \mathbf{1}%
_{A_{n,\delta\varepsilon_{n}}}\left( y_{1}^{n}\right) \frac{\mathfrak{r}%
_{n}\left( y_{1}^{n}\right) }{\mathfrak{s}_{n}(y_{1}^{n})}\mathfrak{s}%
_{n}(y_{1}^{n})dy_{1}^{n}.
\end{equation*}
Since
\begin{equation*}
\mathfrak{R}_{n}\left( A_{n,\delta\varepsilon_{n}}\right) \leq
(1+\delta\varepsilon_{n})\mathfrak{S}_{n}\left(
A_{n,\delta\varepsilon_{n}}\right)
\end{equation*}
it follows that
\begin{equation*}
\lim_{n\rightarrow\infty}\mathfrak{S}_{n}\left(
A_{n,\delta\varepsilon_{n}}\right) =1,
\end{equation*}
which proves the claim.
\end{proof}

Applying this Lemma to the present setting yields%
\[
g_{a_{n}}\left(  Y_{1}\right)  =p_{a_{n}}\left(  Y_{1}\right)  \left(
1+o_{G_{a_{n}}}\left(  1/\sqrt{n}\right)  \right)
\]
as $n\rightarrow\infty.$

This fact entails, as in \cite{Bronia}

\begin{thm}
\label{ThmcvVarTot}Under all the notation and hypotheses above the total
variation norm between $P_{a_{n}}$ and $G_{a_{n}}$ goes to $0$ as
$n\rightarrow\infty.$
\end{thm}

The proof goes as follows

For all $\delta>0$, let
\[
E_{\delta}:=\left\{  y\in\mathbb{R}:\left\vert \frac{p_{a_{n}}\left(
y\right)  -g_{a_{n}}\left(  y\right)  }{g_{a_{n}}\left(  y\right)
}\right\vert <\delta\right\}
\]
which
\begin{equation}
\lim_{n\rightarrow\infty}P_{a_{n}}\left(  E_{\delta}\right)  =\lim
_{n\rightarrow\infty}G_{a_{n}}\left(  E_{\delta}\right)
=1.\label{limPu1,nGu1,n}%
\end{equation}
It holds%
\[
\sup_{C\in\mathcal{B}\left(  \mathbb{R}\right)  }\left\vert P_{a_{n}}\left(
C\cap E_{\delta}\right)  -G_{a_{n}}\left(  C\cap E_{\delta}\right)
\right\vert \leq\delta\sup_{C\in\mathcal{B}\left(  \mathbb{R}\right)  }%
\int_{C\cap E_{\delta}}g_{a_{n}}\left(  y\right)  dy\leq\delta.
\]
By the above result (\ref{limPu1,nGu1,n})
\[
\sup_{C\in\mathcal{B}\left(  \mathbb{R}\right)  }\left\vert P_{a_{n}}\left(
C\cap E_{\delta}\right)  -P_{a_{n}}\left(  C\right)  \right\vert <\eta_{n}%
\]
and
\[
\sup_{C\in\mathcal{B}\left(  \mathbb{R}\right)  }\left\vert G_{a_{n}}\left(
C\cap E_{\delta}\right)  -G_{a_{n}}\left(  C\right)  \right\vert <\eta_{n}%
\]
for some sequence $\eta_{n}\rightarrow0$ ; hence
\[
\sup_{C\in\mathcal{B}\left(  \mathbb{R}\right)  }\left\vert P_{a_{n}}\left(
C\right)  -G_{a_{n}}\left(  C\right)  \right\vert <\delta+2\eta_{n}%
\]
for all positive $\delta,$ which proves the claim.

As a consequence, applying Scheffé's Lemma
\[
\int\left\vert p_{a_{n}}-g_{a_{n}}\right\vert dx\rightarrow0\text{ \ as
}n\rightarrow\infty.
\]
\bigskip

\begin{rem}
This result is to be paralleled with Theorem 1.6 in Diaconis and Freedman
\cite{Diaconis1} and Theorem 2.15 in Dembo and Zeitouni \cite{Dembo} which
provide a rate for this convergence in the LDP range.
\end{rem}

\subsection{\bigskip The asymptotic location of $X$ under the conditioned
distribution}

This section intends to provide some insight on the behaviour of $X_{1}$ under
the condition $\left(  S_{1}^{n}=na_{n}\right)  ;$ this will be extended
further on to the case when $\left(  S_{1}^{n}\geq na_{n}\right)  $ and to be
considered in parallel with similar facts developped in \cite{Bronia} for
larger values of $a_{n}.$

\bigskip It will be seen that conditionally on $\left(  S_{1}^{n}%
=na_{n}\right)  $ the marginal distribution of the sample concentrates around
$a_{n}.$ Let $\mathcal{X}_{t}$ be a r.v. with density $\pi^{a_{n}}$ where
$m(t)=a_{n}$ and $a_{n}$ satisfies (\ref{croissance de a}). Recall that
E$\mathcal{X}_{t}=a_{n}$ $\ and$Var$\mathcal{X}_{t}=s^{2}$. We evaluate the
moment generating function of the normalized variable $\left(  \mathcal{X}%
_{t}-a_{n}\right)  /s$. It holds%
\[
\log E\exp\lambda\left(  \mathcal{X}_{t}-a_{n}\right)  /s=-\lambda
a_{n}/s+\log\Phi\left(  t+\frac{\lambda}{s}\right)  -\log\Phi\left(  t\right)
.
\]
A second order Taylor expansion in the above display yields%
\[
\log E\exp\lambda\left(  \mathcal{X}_{t}-a_{n}\right)  /s=\frac{\lambda^{2}%
}{2}\frac{s^{2}\left(  t+\frac{\theta\lambda}{s}\right)  }{s^{2}}%
\]
where $\theta=\theta(t,\lambda)\in\left(  0,1\right)  .$ It holds

\begin{lem}
Under the above hypotheses and notation, for any compact set $K,$
\[
\lim_{n\rightarrow\infty}\sup_{u\in K}\frac{s^{2}\left(  t+\frac{u}{s}\right)
}{s^{2}}=1.
\]

\end{lem}

Proof: \textbf{Case 1:} if $h(t)\in R_{\beta}$. By Theorem $\ref{order of s}$,
it holds $s^{2}\sim\psi^{\prime}(t)$ with $\psi(t)\sim t^{1/\beta}l_{1}(t)$,
where $l(t)$ is some slowly varying function. And we have also $\psi^{\prime
}(t)=1/h^{^{\prime}}\big(\psi(t)\big)$, hence by $(\ref{3fqg010})$ it follows
\begin{align*}
\frac{1}{s^{2}}  &  \sim h^{^{\prime}}\big(\psi(t)\big)=\psi(t)^{\beta-1}%
l_{0}\big(\psi(t)\big)\big(\beta+\epsilon\big(\psi(t)\big)\big)\\
&  \sim\beta t^{1-1/\beta}l_{1}(t)^{\beta-1}l_{0}\big(\psi(t)\big)=o(t),
\end{align*}
which implies that for any $u\in K$ it holds
\[
\frac{u}{s}=o(t),
\]

\begin{align*}
\frac{s^{2}\left(  t+u/s\right)  }{s^{2}}  &  \sim\frac{\psi^{\prime}%
(t+u/s)}{\psi^{\prime}(t)}=\frac{\psi(t)^{\beta-1}l_{0}\big(\psi
(t)\big)\big(\beta+\epsilon\big(\psi(t)\big)\big)}{\big(\psi
(t+u/s)\big)^{\beta-1}l_{0}\big(\psi(t+u/s)\big)\big(\beta+\epsilon
\big(\psi(t+u/s)\big)\big)}\\
&  \sim\frac{\psi(t)^{\beta-1}}{\psi(t+u/s)^{\beta-1}}\sim\frac{t^{1-1/\beta
}l_{1}(t)^{\beta-1}}{(t+u/s)^{1-1/\beta}l_{1}(t+u/s)^{\beta-1}}\longrightarrow
1.
\end{align*}

\textbf{Case 2:} if $h(t)\in R_{\infty}$. Then we have in this case
$\psi(t)\in\widetilde{R_{0}}$, hence it holds
\[
\frac{1}{st}\sim\frac{1}{t\sqrt{\psi^{\prime}(t)}}=\sqrt{\frac{1}%
{t\psi(t)\epsilon(t)}}\longrightarrow0,
\]
which last step holds from condition $(\ref{3section1030})$. Hence for any
$u\in K$, we get as $n\rightarrow\infty$
\[
\frac{u}{s}=o(t),
\]
thus using the slowly varying propriety of $\psi(t)$ we have
\begin{align*}
\frac{s^{2}\left(  t+u/s\right)  }{s^{2}}  &  \sim\frac{\psi^{\prime}%
(t+u/s)}{\psi^{\prime}(t)}=\frac{\psi(t+u/s)\epsilon(t+u/s)}{t+u/s}\frac
{t}{\psi(t)\epsilon(t)}\\
&  \sim\frac{\epsilon(t+u/s)}{\epsilon(t)}=\frac{\epsilon(t)+O\big(\epsilon
^{\prime}(t)u/s\big)}{\epsilon(t)}\longrightarrow1,
\end{align*}
where we use one Taylor expansion in the second line, and last step holds from
condition $(\ref{3section103})$. This completes the proof.

\bigskip

\bigskip

Applying the above Lemma it follows that the normalized r.v's $\left(
\mathcal{X}_{t}-a_{n}\right)  /s$ converge to a standard normal variable
$N(0,1)$ in distribution, as $n\rightarrow\infty.$ This amount to say that
\[
\mathcal{X}_{t}=a_{n}+sN(0,1)+o_{\Pi^{a_{n}}}(1).
\]
Recall that $\lim_{n\rightarrow\infty}s=0$, which implies that $\mathcal{X}%
_{t}$ concentrates around $a_{n}$ with rate $s.$ Due to Theorem
\ref{ThmcvVarTot} the same holds for $X_{1}$ under $\left(  S_{1}^{n}%
=na_{n}\right)  .$

\subsection{Differences between Gibbs principle under LDP and under extreme
deviation}

It is of interest to confront the present results with the general form of the
Gibbs principle under linear contraints in the LDP range. We recall briefly
and somehow unformally the main classical facts in a simple setting similar as
the one used in this paper.

Let $X_{1},..,X_{n}$ denote $n$ i.i.d. real valued r.v's with distribution $P
$ and density $p$ and let $f:\mathbb{R\rightarrow R}$ be a measurable function
such that $\Phi_{f}(\lambda):=E\exp\lambda f(X_{1})$ is finite for $\lambda$
in a non void neighborhood of $0$ (the so-called Cramer condition). Denote
$m_{f}(\lambda)$ and $s_{f}^{2}(\lambda)$ the first and second derivatives of
$\log\Phi_{f}(\lambda).$ Consider the point set condition $E_{n}:=\left(
\frac{1}{n}\sum_{i=1}^{n}f(X_{i})=0\right)  $ and let $\Omega$ be the set of
all probability measures on $\mathbb{R}$ such that $\int f(x)dQ(x)=0.$

The classical Gibbs conditioning principle writes as follows:

The limiting distribution $P^{\ast}$ of $X_{1}$ conditioned on the family of
events $E_{n}$ exists and is defined as the unique minimizer of the
Kullback-Leibler distance \ between $P$ and $\Omega,$ namely%
\[
P^{\ast}=\arg\min\left\{  K(Q,P),Q\in\Omega\right\}
\]
where
\[
K(Q,P):=\int\log\frac{dQ}{dP}dQ
\]
whenever $Q$ is absolutely continuous w.r.t. $P$, and $K(Q,P)=\infty$
otherwise. \ Also it can be proved that $P^{\ast}$ has a density, which is
defined through
\[
p^{\ast}(x)=\frac{\exp\lambda f(x)}{\Phi_{f}(\lambda)}p(x)
\]
with $\lambda$ the unique solution of the equation $m_{f}(\lambda)=0.$ Take
$f(x)=x-a$ with $a$ fixed to obtain \
\[
p^{\ast}(x)=\pi^{a}(x)
\]
with the current notation of this paper.

Consider now the application of the above result to r.v's $Y_{1},..,Y_{n}$
with $Y_{i}:=\left(  X_{i}\right)  ^{2}$ \ where the $X_{i}^{\prime}s$ are
i.i.d. and are such that the density of the i.i.d. r.v's $Y_{i}^{\prime}s$
satisfy (\ref{densityFunction}) with all the hypothese stated in this paper.
By the Gibbs conditional principle, for \textit{fixed} $a$, conditionally on
$\left(  \sum_{i=1}^{n}Y_{i}=na\right)  $ the generic r.v. $Y_{1}$ has a non
degenerate limit distribution
\[
p_{Y}^{\ast}(y):=\frac{\exp ty}{E\exp tY_{1}}p_{Y}(y)
\]
and the limit density of $X_{1}$ under $\left(  \sum_{i=1}^{n}X_{i}%
^{2}=na\right)  $ is%

\[
p_{X}^{\ast}(y):=\frac{\exp tx^{2}}{E\exp tX_{1}^{2}}p_{X}(y)
\]
whereas, when $a_{n}\rightarrow\infty$ its limit distribution is degenerate
and concentrates around $a_{n}.$ As a consequence the distribution of $X_{1}$
under the condition $\left(  \sum_{i=1}^{n}X_{i}^{2}=na_{n}\right)  $
concentrates sharply at $-\sqrt{a_{n}}$ and $+\sqrt{a_{n}}.$

\section{EDP under exceedance}

The following proposition states the marginally conditional density under
condition $A_{n}=\{S_{n}\geq na_{n}\}$, we denote this density by $p_{A_{n}}$
to differentiate it from $p_{a_{n}}$ which is under condition $\{S_{n}%
=na_{n}\}$. For the purpose of proof, we need the following lemma, based on
Theorem $6.2.1$ of Jensen \cite{Jensen}, to provide one asymptotic estimation
of tail probability $P(S_{n}\geq na_{n})$ and $n$-convolution density
$p(S_{n}/n=u)$ for $u>a_{n}$.

Define
\begin{align}
\label{3virg01000}I(x):=xm^{-1}(x)-\log\Phi\big(m^{-1}(x)\big).
\end{align}

\bigskip

\begin{lem}
\label{JensenLemme} $X_{1},...,X_{n}$ are i.i.d. random variables with density
$p(x)$ defined in $(\ref{densityFunction})$ and $h(x)\in\mathfrak{R}$. Set
$m(t_{n})=a_{n}$. Suppose when $n\rightarrow\infty$, if it holds
\begin{align}
\label{3theorem3cond100}\frac{\psi(t_{n})^{2}}{\sqrt{n}\psi^{\prime}(t_{n}%
)}\longrightarrow0,
\end{align}
then it holds
\begin{align}
\label{3virg0100}P(S_{n}\geq na_{n})=\frac{\exp(-nI(a_{n}))}{\sqrt{2\pi}%
\sqrt{n}t_{n}s(t_{n})}\Big(1+o\big(\frac{1}{\sqrt{n}}\big)\Big).
\end{align}
Let further
\[
H_{n}(u):=\frac{\sqrt{n}\exp\big(-nI(u)\big)}{\sqrt{2\pi}s(t_{u})}%
\]
It then holds
\begin{align}
\label{3virg01001}\sup_{u>a_{n}}\frac{p(S_{n}/n=u)}{H_{n}(u)}=1+o\left(
1/\sqrt{n}\right)  .
\end{align}

\end{lem}

Proof: For the density $p(x)$ defined in $(\ref{densityFunction})$, we show
$g(x)$ is convex when $x$ is large enough. If $h(x)\in R_{\beta}$, it holds
for $x$ large enough
\begin{align}
\label{3virg01}g^{^{\prime\prime}}(x)=h^{^{\prime}}(x)=\frac{h(x)}%
{x}\big(\beta+\epsilon(x)\big)>0.
\end{align}
If $h(x)\in R_{\infty}$, its reciprocal function $\psi(x)\in\widetilde{R_{0}}%
$. Set $x=\psi(u)$, hence we have for $x$ large enough
\begin{align}
\label{3virg02}g^{^{\prime\prime}}(x)=h^{^{\prime}}(x)=\frac{1}{\psi
^{^{\prime}}(u)}=\frac{u}{\psi(u)\epsilon(u)}>0,
\end{align}
where the inequality holds since $\epsilon(u)>0$ under condition
$(\ref{3section1030})$ when $u$ is large enough. $(\ref{3virg01})$ and
$(\ref{3virg02})$ imply that $g(x)$ is convex for $x$ large enough.

Therefore, the density $p(x)$ with $h(x)\in\mathfrak{R}$ satisfies the
conditions of Jensen's Theorem 6.2.1 (\cite{Jensen}). Denote by $p_{n}$ the
density of $\bar{X}=(X_{1}+...+X_{n})/n$. We obtain with the third order's
Edgeworth expansion from formula $(2.2.6)$ of $(\cite{Jensen})$
\begin{align}
\label{3theoreom1010}P(S_{n}\ge n a_{n})=\frac{\Phi(t_{n})^{n}\exp
(-nt_{n}a_{n})}{\sqrt{n}t_{n} s(t_{n})}\Big(B_{0}(\lambda_{n})+O\big(\frac
{\mu_{3}(t_{n})}{6\sqrt{n}s^{3}(t_{n})}B_{3}(\lambda_{n})\big)\Big),
\end{align}
where $\lambda_{n}=\sqrt{n}t_{n}s(t_{n})$, $B_{0}(\lambda_{n})$ and
$B_{3}(\lambda_{n})$ are defined by
\begin{align*}
B_{0}(\lambda_{n})=\frac{1}{\sqrt{2\pi}}\Big(1-\frac{1}{\lambda_{n}^{2}%
}+o(\frac{1}{\lambda_{n}^{2}})\Big),\qquad B_{3}(\lambda_{n})\sim-\frac
{3}{\sqrt{2\pi}\lambda_{n}}.
\end{align*}
We show, under condition $(\ref{3theorem3cond100})$, it holds as
$a_{n}\rightarrow\infty$
\begin{align}
\label{3theoreom11}\frac{1}{\lambda_{n}^{2}}=o\big(\frac{1}{n}\big).
\end{align}
Since $n/\lambda_{n}^{2}=1/(t_{n}^{2}s^{2}(t_{n}))$, $(\ref{3theoreom11})$ is
equivalent to show
\begin{align}
\label{3theoreom112}t_{n}^{2}s^{2}(t_{n})\longrightarrow\infty.
\end{align}
By Theorem $\ref{order of s}$, $m(t_{n})\sim\psi(t_{n})$ and $s^{2}(t_{n}%
)\sim\psi^{\prime}(t_{n})$, combined with $(\ref{3theoreom110})$, it holds
$t_{n}\sim h(a_{n})$.

If $h\in R_{\beta}$, notice that it holds
\begin{align*}
\psi^{\prime}(t_{n})=\frac{1}{h^{\prime}(\psi(t_{n}))}=\frac{\psi(t_{n}%
)}{h\big(\psi(t_{n})\big)\big(\beta+\epsilon(\psi(t_{n}))\big)}\sim\frac
{a_{n}}{h(a_{n})\big(\beta+\epsilon(\psi(t_{n}))\big)},
\end{align*}
hence we have
\begin{align}
\label{3lh01}t_{n}^{2}s^{2}(t_{n})\sim h(a_{n})^{2} \frac{a_{n}}%
{h(a_{n})\big(\beta+\epsilon(\psi(t_{n}))\big)}= \frac{a_{n}h(a_{n})}%
{\beta+\epsilon(\psi(t_{n}))}\longrightarrow\infty.
\end{align}

If $h\in R_{\infty}$, then $\psi(t_{n})\in\widetilde{R_{0}}$, thus it follows
\begin{align}
\label{3lh010}t_{n}^{2}s^{2}(t_{n})\sim t_{n}^{2} \frac{\psi(t_{n}%
)\epsilon(t_{n})}{t_{n}}=t_{n}\psi(t_{n})\epsilon(t_{n})\longrightarrow\infty,
\end{align}
where last step holds from condition $(\ref{3section1030})$. We have showed
$(\ref{3theoreom11})$ , therefore it holds
\begin{align*}
B_{0}(\lambda_{n})=\frac{1}{\sqrt{2\pi}}\Big(1+o(\frac{1}{n})\Big).
\end{align*}

By $(\ref{3theoreom112})$, $\lambda_{n}$ goes to $\infty$ as $a_{n}%
\rightarrow\infty$, which implies further $B_{3}(\lambda_{n})\rightarrow0$. On
the other hand, by $(\ref{3cor1})$ it holds $\mu_{3}/s^{3}\rightarrow0$. Hence
we obtain from $(\ref{3theoreom1010})$
\[
P(S_{n}\geq na_{n})=\frac{\Phi(t_{n})^{n}\exp(-nt_{n}a_{n})}{\sqrt{2\pi
n}t_{n}s(t_{n})}\Big(1+o\big(\frac{1}{\sqrt{n}}\big)\Big),
\]
which together with $(\ref{3virg01000})$ gives $(\ref{3virg0100})$.

By Jensen's Theorem $6.2.1$ (\cite{Jensen}) and formula $(2.2.4)$
in\cite{Jensen} it follows that
\[
p(S_{n}=na_{n})=\frac{\sqrt{n}\Phi(t_{n})^{n}\exp(-nt_{n}a_{n})}{\sqrt{2\pi
}s(t_{n})}\Big(1+o\big(\frac{1}{\sqrt{n}}\big)\Big),
\]

\bigskip which, together with $(\ref{3virg01000})$, gives $(\ref{3virg01001})$.

\begin{prop}
\label{3theorem5} \label{main theorem} $X_{1},...,X_{n}$ are i.i.d. random
variables with density $p(x)$ defined in $(\ref{densityFunction})$ and
$h(x)\in\mathfrak{R}$. Suppose when $n\rightarrow\infty$, if it holds
\begin{align}
\label{3theorem3cond1000}\frac{\psi(t_{n})^{2}}{\sqrt{ n}\psi^{\prime}(t_{n})}
\longrightarrow0,
\end{align}
and
\begin{align}
\label{edgeworth condition01}\eta_{n} \rightarrow0, \qquad\qquad\frac{\log
n}{nh(a_{n}) \eta_{n}} \rightarrow\infty,
\end{align}
then
\begin{align*}
p_{A_{n}}(y_{1})=p(X_{1}=y_{1} \vert S_{n} \ge n a_{n})=g_{A_{n}}%
(y_{1})\Big(1+o\big(\frac{1}{\sqrt{n}}\big)\Big),
\end{align*}
where $g_{A_{n}}(y_{1})=nt_{n} s(t_{n})e^{nI(a_{n})}\int_{a_{n}}^{a_{n}%
+\eta_{n}} g_{\tau}(y_{1})\exp\big(-nI(\tau)-\log s(t_{\tau})\big) d\tau$,
$g_{\tau}(y_{1})$ is defined as $g_{a_{n}}(y_{1})$ in Theorem
$(\ref{point conditional density})$ on replacing $a_{n}$ by $\tau$.
\end{prop}

Proof: We can denote $p_{A_{n}}(y_{1})$ by the integration of $p_{a_{n}}%
(y_{1})$
\begin{align*}
p_{A_{n}}(y_{1})  &  =\int_{a_{n}}^{\infty}p(X_{1}=y_{1}\vert S_{n}=n
\tau)p(S_{n}=n\tau\vert S_{n} \ge n a_{n}) d\tau\\
&  =p(X_{1}=y_{1})\frac{\int_{a_{n}}^{\infty}p(S_{2}^{n}=n \tau-y_{1})d \tau
}{p(S_{n} \ge n a_{n})}\\
&  =\frac{p(X_{1}=y_{1})}{p(S_{n} \ge n a_{n})} P_{1}\Big(1+\frac{P_{2}}%
{P_{1}}\Big),
\end{align*}
where the second equality is obtained by Bayes formula, and $P_{1}=\int
_{a_{n}}^{a_{n}+\eta_{n}} p(S_{2}^{n}=n \tau-y_{1})d \tau$, $P_{2}=\int
_{a_{n}+\eta_{n}}^{\infty}p(S_{2}^{n}=n \tau-y_{1})d \tau$, $S_{2}^{n}%
=X_{2}+...+X_{n}$. In fact $P_{2}$ is one infinitely small term with respect
to $P_{1}$, which is proved below. Further we have
\begin{align*}
&  P_{2}=\frac{1}{n}P\Big( S_{2}^{n}\ge{n(a_{n}+\eta)-y_{1}}\Big)=\frac{1}%
{n}P\Big( S_{2}^{n}\ge(n-1)c_{n}\Big),\\
&  P_{1}+P_{2}=\frac{1}{n}P\Big( S_{2}^{n}\ge{na_{n}-y_{1}}\Big)=\frac{1}%
{n}P\Big( S_{2}^{n}\ge(n-1)d_{n}\Big),
\end{align*}
where $c_{n}=\big(n(a_{n}+\eta_{n})-y_{1}\big)/(n-1)$ and $d_{n}=(na_{n}%
-y_{1})/(n-1)$. Denote $t_{c_{n}}=m^{-1}(c_{n})$ and $t_{d_{n}}=m^{-1}(d_{n}%
)$. Using Lemma $(\ref{JensenLemme})$, it holds
\begin{align}
\label{3the51}\frac{P_{2}}{P_{1}+P_{2}}=\Big(1+o\big(\frac{1}{\sqrt{n}%
}\big)\Big)\frac{t_{d_{n}}s(t_{d_{n}})}{t_{c_{n}}s(t_{c_{n}})}\exp
\Big(-(n-1)\big(I(c_{n})-I(d_{n})\big)\Big),
\end{align}
Using the convexity of the function $I$, it holds
\begin{align*}
\exp\Big(-(n-1)\big(I(c_{n})-I(d_{n})\big)\Big)  &  \le\exp\Big(-(n-1)(c_{n}%
-d_{n})m^{-1}(d_{n})\big)\Big)\\
&  =\exp\big(-n\eta_{n} m^{-1}(d_{n})\big)
\end{align*}
Consider $u\rightarrow m^{-1}(u)$ is increasing, since $d_{n} \le a_{n}$ as
$a_{n}\rightarrow\infty$, it holds $m^{-1}(d_{n})\ge m^{-1}(a_{n})$, hence we
get
\begin{align}
\label{3lh002}\exp\Big(-(n-1)\big(I(c_{n})-I(d_{n})\big)\Big)  &  \le
\exp\big(-n\eta_{n} m^{-1}(a_{n})\big).
\end{align}
Using Theorem $\ref{order of s}$, we have $m^{-1}(a_{n})\sim\psi^{-1}%
(a_{n})=h(a_{n})$, thus under condition $(\ref{edgeworth condition01})$ it
holds as $a_{n}\rightarrow\infty$
\begin{align*}
\exp\Big(-(n-1)\big(I(c_{n})-I(d_{n})\big)\Big) \longrightarrow0.
\end{align*}

Then we show it holds
\begin{align}
\label{3lh502}\frac{t_{d_{n}}s(t_{d_{n}})}{t_{c_{n}}s(t_{c_{n}})}%
\longrightarrow1.
\end{align}
By definition, $c_{n}/d_{n}\rightarrow1$ as $a_{n}\rightarrow\infty$. if $h\in
R_{\beta}$, by $(\ref{3lh01})$, it holds
\begin{align}
\label{3lh02}\Big(\frac{t_{d_{n}}s(t_{d_{n}})}{t_{c_{n}}s(t_{c_{n}})}\Big)^{2}
\sim\Big(\frac{d_{n} h(d_{n})}{\beta+\epsilon\big(\psi(d_{n})\big)}\Big)^{2}
\Big(\frac{\beta+\epsilon\big(\psi(c_{n})\big)}{c_{n} h(c_{n})}\Big)^{2}%
\sim\Big(\frac{h(d_{n})}{h(c_{n})}\Big)^{2} \longrightarrow1.
\end{align}
If $h\in R_{\infty}$, notice the function $t\rightarrow t\psi(t)\epsilon(t)$
is increasing and continuous as $t$ large enough. By $(\ref{3lh010})$, it
holds
\begin{align}
t^{2}s^{2}(t)\sim t\psi(t)\epsilon(t),
\end{align}
consider $d_{n} \rightarrow c_{n}$ as $n\rightarrow\infty$, hence we have
\begin{align}
\Big(\frac{t_{d_{n}}s(t_{d_{n}})}{t_{c_{n}}s(t_{c_{n}})}\Big)^{2} \sim
\frac{d_{n}\psi(d_{n})\epsilon(d_{n})}{c_{n}\psi(c_{n})\epsilon(c_{n})}
\longrightarrow1.
\end{align}
Using $(\ref{3the51})$, $(\ref{3lh002})$ and $(\ref{3lh502})$, we obtain
\begin{align*}
\frac{P_{2}}{P_{1}+P_{2}} \le2 \exp\big(-n m^{-1}(a_{n})\eta_{n}\big) ,
\end{align*}
which, together with condition $(\ref{edgeworth condition01})$, it holds
\begin{align*}
\frac{P_{2}}{P_{1}} =o\big(\frac{1}{\sqrt{n}}\big).
\end{align*}

Therefore we can approximate $p_{A_{n}}(y_{1})$ by
\begin{align}
\label{3lfl0510}p_{A_{n}}(y_{1})=\Big(1+o\big(\frac{1}{\sqrt{n}}%
\big)\Big)\int_{a_{n}}^{a_{n}+\eta_{n}}p(X_{1}=y_{1}|S_{n}=n\tau)p(S_{n}%
=n\tau|S_{n}\geq na_{n})d\tau.
\end{align}
According to Lemma $\ref{JensenLemme}$, it follows when $\tau\in\lbrack
a_{n},a_{n}+\eta_{n}]$
\begin{align}
\label{3lfl010}p(S_{n}=n\tau|S_{n}\geq na_{n})=\Big(1+o\big(\frac{1}{\sqrt{n}%
}\big)\Big)\frac{nm^{-1}(a_{n})s(t_{n})}{s(t_{\tau})}\exp\big(-n(I(\tau
)-I(a_{n}))\big),
\end{align}
where $m(t_{n})=a_{n}$, $m(t_{\tau})=\tau$. Inserting $(\ref{3lfl0510})$ into
$(\ref{3lfl010})$, we obtain
\[
p_{A_{n}}(y_{1})=\Big(1+o\big(\frac{1}{\sqrt{n}}\big)\Big)nt_{n}%
s(t_{n})e^{nI(a_{n})}\int_{a_{n}}^{a_{n}+\eta_{n}}g_{\tau}(y_{1}%
)\exp\big(-nI(\tau)-\log s(t_{\tau})\big)d\tau,
\]
this completes the proof.

\bigskip

\section{Appendix}

For density functions $p(x)$ defined in $(\ref{densityFunction})$ satisfying
also $h(x)\in\mathfrak{R}$, denote by $\psi(x)$ the reciprocal function of
$h(x)$ and $\sigma^{2}(v)=\big(h^{\prime}(v)\big)^{-1}$, $v\in\mathbb{R}_{+}$.
For brevity, we write $\hat{x},\sigma,l$ instead of $\hat{x}(t),\sigma
\big(\psi(t)\big),l(t)$.

When $t$ is given, $K(x,t)$ attain its maximum at $\hat{x}=\psi(t)$. The
fourth order Taylor expansion of $K(x,t)$ on $x\in[\hat{x}-\sigma l,\hat
{x}+\sigma l]$ yields
\begin{align}
\label{3abeltheorem001}K(x,t)=K(\hat{x},t)-\frac{1}{2}h^{\prime}(\hat
{x})\big(x-\hat{x}\big)^{2} -\frac{1}{6}h^{\prime\prime}(\hat{x}%
)\big(x-\hat{x}\big)^{3}+\epsilon(x,t),
\end{align}
with some $\theta\in[0,1]$
\begin{align}
\label{3abel01}\epsilon(x,t)=-\frac{1}{24}h^{^{\prime\prime\prime}}%
\big(\hat{x}+\theta(x-\hat{x})\big)(x-\hat{x})^{4}.
\end{align}

\begin{lem}
\label{3lemma00} For $p(x)$ in $(\ref{densityFunction})$, $h(x)\in
\mathfrak{R}$, it holds when $t\rightarrow\infty$,
\begin{align}
\label{3section201}\frac{|\log\sigma\big(\psi(t)\big)|}{\int_{1}^{t}
\psi(u)du}\longrightarrow0.
\end{align}

\end{lem}

Proof: If $h(x)\in R_{\beta}$, by Theorem $(1.5.12)$ of $\cite{Bingham}$,
there exists some slowly varying function such that it holds $\psi(x)\sim
x^{1/\beta}l_{1}(x) $. Hence it holds as $t\rightarrow\infty$(see
\cite{Feller}, Chapter 8)
\begin{align}
\label{3section202}\int_{1}^{t} \psi(u)du \sim t^{1+\frac{1}{\alpha}}
l_{1}(t).
\end{align}
On the other hand, $h^{\prime}(x)= x^{\beta-1}l(x)\big(\beta+\epsilon
(x)\big)$, thus we have as $x\rightarrow\infty$
\begin{align*}
|\log\sigma(x)|  &  =\big|\log\big(h^{\prime}(x)\big)^{-\frac{1}{2}%
}\big|=\Big|\frac{1}{2}\big((\beta-1)\log x+\log l(x)+\log(\beta
+\epsilon(x))\big)\Big|\\
&  \le\frac{1}{2}(\beta+1)\log x,
\end{align*}
set $x=\psi(t)$, then when $t\rightarrow\infty$, it holds $x<2 t^{1/\beta
}l_{1}(t)<t^{1/\beta+1}$, hence we have
\begin{align*}
|\log\sigma\big(\psi(t)\big)|<\frac{(\beta+1)^{2}}{2\beta} \log t,
\end{align*}
which, together with $(\ref{3section202})$, yields $(\ref{3section201})$.

If $h(x)\in R_{\infty}$, then by definition $\psi(x)\in\widetilde{R_{0}}$ is
slowly varying as $x\rightarrow\infty$. Hence it holds as $t\rightarrow\infty
$(see \cite{Feller}, Chapter 8)
\begin{align}
\label{3section203}\int_{1}^{t} \psi(u)du \sim t\psi(t).
\end{align}
And now we have $h^{\prime}(x)=1/\psi^{\prime}(v)$ with $x=\psi(v)$. Therefore
it follows
\begin{align*}
|\log\sigma(x)|=\big|\log\big(h^{\prime}(x)\big)^{-\frac{1}{2}}\big|=\frac
{1}{2} |\log\psi^{\prime}(v)|,
\end{align*}
Set $x=\psi(t)$, then $v=t$, consider $\psi(t)\in\widetilde{R_{0}}$, thus we
have
\begin{align}
\label{3section204}|\log\sigma\big(\psi(t)\big)|  &  =\frac{1}{2} |\log
\psi^{\prime}(t)|=\frac{1}{2} \Big|\log\Big(\psi(t)\frac{\epsilon(t)}%
{t}\Big)\Big|\nonumber\\
&  =\frac{1}{2} \big|\log\psi(t)+\log\epsilon(t)-\log{t}\big|\nonumber\\
&  \le\log t+\frac{1}{2} |\log\epsilon(t)|\le2\log t,
\end{align}
where last inequality follows from $(\ref{3section103})$. $(\ref{3section203}%
)$ and $(\ref{3section204})$ imply $(\ref{3section201})$. This completes the proof.

\begin{lem}
\label{3lemma0} For $p(x)$ in $(\ref{densityFunction})$, $h\in\mathfrak{R}$,
then for any varying slowly function $l(t)\rightarrow\infty$ as $t\rightarrow
\infty$, it holds
\begin{align}
\label{3section201}\sup_{|x|\le\sigma l}{h^{\prime\prime\prime}(\hat{x}%
+x)}\sigma^{4} l^{4}\longrightarrow0 \qquad as \quad t\rightarrow\infty.
\end{align}

\end{lem}

Proof: \textbf{Case 1:} $h\in R_{\beta}$. We have $h(x)=x^{\beta}l_{0}(x) ,
l_{0}(x)\in R_{0},\beta>0$. Hence it holds
\begin{align}
\label{mu3 10101}h^{^{\prime\prime}}(x)=\beta(\beta-1)x^{\beta-2}%
l_{0}(x)+2\beta x^{\beta-1}l_{0}^{^{\prime}}(x)+x^{\beta}l_{0}^{^{\prime
\prime}}(x).
\end{align}
and
\begin{align}
\label{mu3 1010}h^{^{\prime\prime\prime}}(x)=\beta(\beta-1)(\beta
-2)x^{\beta-3}l_{0}(x)+3\beta(\beta-1)x^{\beta-2}l_{0}^{^{\prime}}(x)+3\beta
x^{\beta-1}l_{0}^{^{\prime\prime}}(x)+x^{\beta}l_{0}^{^{\prime\prime\prime}%
}(x).
\end{align}
Consider $l(x)\in R_{0}$, it is easy to obtain
\begin{align}
\label{mu3 10102}l_{0}^{^{\prime}}(x)=\frac{l_{0}(x)}{x}\epsilon(x),\qquad
l_{0}^{^{\prime\prime}}(x)=\frac{l_{0}(x)}{x^{2}}\big(\epsilon^{2}%
(x)+x\epsilon^{\prime}(x)-\epsilon(x)\big),
\end{align}
and
\begin{align*}
l_{0}^{^{\prime\prime\prime}}(x)=\frac{l_{0}(x)}{x^{3}}\big(\epsilon
^{3}(x)+3x\epsilon^{\prime}(x)\epsilon(x)-3\epsilon^{2}(x) -2x\epsilon
^{^{\prime}}(x)+2\epsilon(x)+x^{2}\epsilon^{^{\prime\prime}}(x)\big).
\end{align*}
Under condition $(\ref{3section104})$, there exists some positive constant $Q
$ such that it holds
\begin{align*}
|l_{0}^{^{\prime\prime}}(x)|\le Q\frac{l_{0}(x)}{x^{2}},\qquad|l_{0}%
^{^{\prime\prime\prime}}(x)|\le Q\frac{l_{0}(x)}{x^{3}},
\end{align*}
which, together with $(\ref{mu3 1010})$, yields with some positive constant
$Q_{1}$
\begin{align}
\label{mu3 1011}|h^{^{\prime\prime\prime}}(x)|\le Q_{1} \frac{h(x)}{x^{3}}.
\end{align}
By definition, we have $\sigma^{2}(x)=1/h^{^{\prime}}(x)=x/\big(h(x)(\beta
+\epsilon(x))\big)$, thus it follows
\begin{align}
\label{mu3 1013}\sigma^{2}=\sigma^{2}(\hat{x})=\frac{\hat{x}}{h(\hat{x}%
)(\beta+\epsilon(\hat{x}))}=\frac{\psi(t)}{t(\beta+\epsilon(\psi(t)))}%
=\frac{\psi(t)}{\beta t}\big(1+o(1)\big),
\end{align}
this implies $\sigma l=o(\psi(t))=o(\hat{x})$. Thus we get with
$(\ref{mu3 1011})$
\begin{align}
\label{mu3 1014}\sup_{|x|\le\sigma l}|h^{^{\prime\prime\prime}}(\hat{x}%
+x)|\le\sup_{|x|\le\sigma l} Q_{1} \frac{h(\hat{x}+x)}{(\hat{x}+x)^{3}}\le
Q_{2} \frac{t}{\psi^{3}(t)},
\end{align}
where $Q_{2}$ is some positive constant. Combined with $(\ref{mu3 1013})$, we
obtain
\begin{align*}
\sup_{|x|\le\sigma l}|h^{^{\prime\prime\prime}}(\hat{x}+x)|\sigma^{4} l^{4}
\le Q_{2} \frac{t}{\psi^{3}(t)}\sigma^{4}l^{4}=\frac{Q_{2}l^{4}}{\beta^{2} t
\psi(t)}\longrightarrow0,
\end{align*}
as sought.

\textbf{Case 2:} $h\in R_{\infty}$. Since $\hat{x}=\psi(t)$, we have
$h(\hat{x})=t$. Thus it holds
\begin{align}
\label{mu3 0010}h^{\prime}(\hat{x})=\frac{1}{\psi^{\prime}(t)}\qquad and \quad
h^{\prime\prime}(\hat{x})=-\frac{\psi^{\prime\prime}(t)}{\big(\psi^{\prime
}(t)\big)^{3}},
\end{align}
further we get
\begin{align}
\label{mu3 001}h^{\prime\prime\prime}(\hat{x})=-\frac{\psi^{^{\prime
\prime\prime}}(t)\psi^{^{\prime}}(t)-3\big(\psi^{^{\prime\prime}}(t)\big)^{2}%
}{\big(\psi^{\prime}(t)\big)^{4}}.
\end{align}
Notice if $h(\hat{x})\in R_{\infty}$, then $\psi(t) \in\widetilde{R_{0}}$.
Therefore we obtain
\begin{align}
\label{mu3 002}\psi^{^{\prime}}(t)=\frac{\psi(t)}{t}\epsilon(t),
\end{align}
and
\begin{align}
\label{mu3 003}\psi^{^{\prime\prime}}(t)  &  =-\frac{\psi(t)}{t^{2}}%
\epsilon(t)\Big(1-\epsilon(t)-\frac{t\epsilon^{^{\prime}}(t)}{\epsilon
(t)}\Big)\nonumber\\
&  =-\frac{\psi(t)}{t^{2}}\epsilon(t)\big(1+o(1)\big) \qquad as \quad
t\rightarrow\infty,
\end{align}
where last equality holds from $(\ref{3section103})$. Using
$(\ref{3section103})$ once again, we have also $\psi^{^{\prime\prime\prime}%
}(t)$
\begin{align}
\label{mu3 004}\psi^{^{\prime\prime\prime}}(t)  &  =\frac{\psi(t)}{t^{3}%
}\epsilon(t)\Big(2+\epsilon^{2}(t)+3t\epsilon^{^{\prime}}(t)-3\epsilon
(t)-\frac{2t\epsilon^{^{\prime}}(t)}{\epsilon(t)}+\frac{t^{2}\epsilon
^{^{\prime\prime}}(t)}{\epsilon(t)}\Big)\nonumber\\
&  =\frac{\psi(t)}{t^{3}}\epsilon(t)\big(2+o(1)\big) \qquad as \quad
t\rightarrow\infty.
\end{align}
Put $(\ref{mu3 002})$ $(\ref{mu3 003})$ and $(\ref{mu3 004})$ into
$(\ref{mu3 001})$ we get
\begin{align*}
h^{^{\prime\prime\prime}}(\hat{x})=\frac{1}{\psi^{2}(t)\epsilon^{2}%
(t)}\big(1+o(1)\big)
\end{align*}
Thus by $(\ref{3section1030})$ it holds as $t\rightarrow\infty$
\begin{align}
\label{mu3 005}\sup_{|v|\le t/4}h^{\prime\prime\prime}\big(\psi(t+v)\big)  &
=\sup_{|v|\le t/4}\frac{1}{\psi^{2}(t+v)\epsilon^{2}(t+v)}%
\big(1+o(1)\big)\nonumber\\
&  \le\sup_{|v|\le t/4}\frac{2\sqrt{t+v}}{\psi^{2}(t+v)}\le\frac{3\sqrt{t}%
}{\psi^{2}(t)},
\end{align}
where last inequality holds from the slowly varying propriety: $\psi
(t+v)\sim\psi(t)$. Using $\sigma=\big(h^{^{\prime}}(\hat{x})\big)^{-1/2}$, it
holds
\begin{align*}
\sup_{|v|\le t/4}{h^{^{\prime\prime\prime}}\big(\psi(t+v)\big)}\sigma^{4}
\le\frac{3\sqrt{t}}{\psi^{2}(t)}\frac{1}{({h^{^{\prime}}(\hat{x})})^{2}}%
=\frac{3\sqrt{t}}{\psi^{2}(t)} \frac{\psi^{2}(t)\epsilon^{2}(t)}{t^{2}}%
=\frac{3\epsilon^{2}(t)}{t^{3/2}}\longrightarrow0,
\end{align*}
where $\epsilon(t)\rightarrow0$ and $\psi(t)$ varies slowly. Hence for any
slowly varying function $l(t)\rightarrow\infty$ it holds as $t\rightarrow
\infty$
\begin{align*}
\sup_{|v|\le t/4}{h^{^{\prime\prime\prime}}\big(\psi(t+v)\big)} \sigma^{4}
l^{4} \longrightarrow0.
\end{align*}
Consider $\psi(t)\in\widetilde{R_{0}}$, thus $\psi(t)$ is increasing, we have
the relation
\begin{align*}
\sup_{|v|\le t/4}{h^{^{\prime\prime\prime}}\big(\psi(t+v)\big)}=\sup
_{|\zeta|\le[\zeta_{1},\zeta_{2}]}h^{^{\prime\prime\prime}}(\hat{x}+\zeta),
\end{align*}
where
\begin{align*}
\zeta_{1}=\psi(3t/4)-\hat{x}, \qquad\zeta_{2}=\psi(5t/4)-\hat{x}.
\end{align*}
Hence we have showed
\begin{align*}
\sup_{|\zeta|\le[\zeta_{1},\zeta_{2}]}{h^{^{\prime\prime\prime}}(\hat{x}%
+\zeta)} \sigma^{4} l^{4} \longrightarrow0.
\end{align*}
For completing the proof, it remains to show
\begin{align}
\label{mu3 0060}\sigma l \le\min(|\zeta_{1}|,\zeta_{2}) \qquad as \quad
t\rightarrow\infty.
\end{align}
Perform first order Taylor expansion of $\psi(3t/4)$ at $t$, for some
$\alpha\in[0,1]$, it holds
\begin{align*}
\zeta_{1}  &  =\psi(3t/4)-\hat{x}=\psi(3t/4)-\psi(t)=-\psi^{^{\prime}%
}\big(t-\alpha{t}/{4}\big) \frac{t}{4}=-\frac{\psi\big(t-\alpha t/{4}%
\big)}{4-\alpha}\epsilon\big(t-\alpha t/{4}\big),
\end{align*}
thus using $(\ref{3section1030})$ and slowly varying propriety of $\psi(t)$ we
get as $t\rightarrow\infty$
\begin{align}
\label{mu3 007}|\zeta_{1}|\ge\frac{\psi\big(t-\alpha t/{4}\big)}{4}%
\epsilon\big(t-\alpha t/{4}\big)\ge\frac{\psi(t)}{5}\epsilon\big(t-\alpha
t/{4}\big)\ge\frac{\psi(t)}{5t^{1/4}}.
\end{align}
On the other hand, we have $\sigma=\big(h^{^{\prime}}(\hat{x})\big)^{-1/2}%
=\big(\psi(t)\epsilon(t)/t\big)^{1/2}$, which, together with $(\ref{mu3 007}%
)$, yields
\begin{align*}
\frac{\sigma}{|\zeta_{1}|}\le5\sqrt{\frac{\epsilon(t)}{\psi(t)\sqrt{t}}%
}\longrightarrow0 \qquad as \quad t\rightarrow\infty,
\end{align*}
which implies for any slowly varying function $l(t)$ it holds $\sigma
l=o(|\zeta_{1}|)$. By the same way, it is easy to show $\sigma l=o(\zeta_{2}%
)$. Hence $(\ref{mu3 0060})$ holds, as sought.

\begin{lem}
\label{3lemma01} For $p(x)$ in $(\ref{densityFunction})$, $h\in\mathfrak{R}$,
then for any varying slowly function $l(t)\rightarrow\infty$ as $t\rightarrow
\infty$, it holds
\begin{align}
\label{3section201}\sup_{|x|\le\sigma l}\frac{h^{\prime\prime\prime}(\hat
{x}+x)}{h^{\prime\prime}(\hat{x})}\sigma l\longrightarrow0 \qquad as \quad
t\rightarrow\infty.
\end{align}
and
\begin{align}
\label{3section2010}h^{^{\prime\prime}}(\hat{x})\sigma^{3} l \longrightarrow0.
\end{align}

\end{lem}

Proof: \textbf{Case 1:} Using $(\ref{mu3 10101})$ and $(\ref{mu3 10102})$, we
get $h^{^{\prime\prime}}(x)=\big(\beta(\beta-1)+o(1)\big)x^{\beta-2} l_{0}(x)$
as $x\rightarrow\infty$, where $l_{0}(x)\in R_{0}$. Hence it holds
\begin{align}
\label{3lem230}h^{^{\prime\prime}}(\hat{x})=\big(\beta(\beta-1)+o(1)\big)\psi
(t)^{\beta-2} l_{0}(\psi(t)),
\end{align}
which, together with $(\ref{mu3 1013} )$ and $(\ref{mu3 1014} )$, yields with
some positive constant $Q_{3}$
\begin{align*}
\sup_{|x|\le\sigma l}\Big|\frac{h^{\prime\prime\prime}(\hat{x}+x)}%
{h^{\prime\prime}(\hat{x})}\sigma l\Big|\le Q_{3} \frac{t}{\psi^{3}(t)}
\frac{1}{\psi(t)^{\beta-2} l_{0}(\psi(t))}\sqrt{\frac{\psi(t)}{\beta t}}l
=\frac{Q_{3}}{\sqrt{\beta}} \frac{\sqrt{t}}{\psi(t)^{\beta+1/2}l_{0}(\psi
(t))}l.
\end{align*}
Notice $\psi(t)\sim t^{1/\beta}l_{1}(t)$ for some slowly varying function
$l_{1}(t)$, then it holds $\sqrt{t}l=o\big(\psi(t)^{\beta+1/2}\big)$. Hence we
get $(\ref{3section201})$.

From $(\ref{mu3 1013} )$ and $(\ref{3lem230})$, we obtain as $t\rightarrow
\infty$
\begin{align}
\label{3hfe1}h^{^{\prime\prime}}(\hat{x})\sigma^{3} l  &  =\big(\beta
(\beta-1)+o(1)\big)\psi(t)^{\beta-2} l_{0}(\psi(t))\Big(\frac{\psi(t)}{\beta
t}\Big)^{3/2}l\nonumber\\
&  =\big(\beta(\beta-1)+o(1)\big) \frac{\psi(t)^{\beta-1/2}}{\beta^{3/2}
t^{3/2}}l_{0}(\psi(t))l \le\frac{1}{\sqrt{t}},
\end{align}
where last inequality holds since $\psi(t)^{\beta-1/2}/t^{3/2}\sim
l_{1}(t)^{\beta-1/2}/t^{1/2+1/2\beta}$ as $t\rightarrow\infty$. This implies
$(\ref{3section2010})$ holds.

\textbf{Case 2:} Using $(\ref{mu3 0010})$ and $(\ref{mu3 003})$ we obtain
\begin{align}
\label{mu3 006}h^{\prime\prime}(\hat{x})=-\frac{\psi^{\prime\prime}%
(t)}{\big(\psi^{\prime}(t)\big)^{3}} =\frac{t}{\psi^{2}(t)\epsilon^{2}%
(t)}\big(1+o(1)\big).
\end{align}
Combine $(\ref{mu3 005})$ and $(\ref{mu3 006})$, using $\sigma
=\big(h^{^{\prime}}(\hat{x})\big)^{-1/2}$, we have as $t\rightarrow\infty$
\begin{align*}
\sup_{|v|\le t/4}\frac{h^{^{\prime\prime\prime}}\big(\psi(t+v)\big)}
{h^{^{\prime\prime}}(\hat{x})}\sigma\le\frac{4\epsilon^{2}(t)}{\sqrt{t}}%
\frac{1}{\sqrt{h^{^{\prime}}(\hat{x})}}=\frac{4\epsilon(t)^{5/2}\sqrt{\psi
(t)}}{t}\rightarrow0,
\end{align*}
where $\epsilon(t)\rightarrow0$ and $\psi(t)$ varies slowly. Hence for
arbitrarily slowly varying function $l(t)$ it holds as $t\rightarrow\infty$
\begin{align*}
\sup_{|v|\le t/4}\frac{h^{^{\prime\prime\prime}}\big(\psi(t+v)\big)}
{h^{^{\prime\prime}}(\hat{x})}\sigma l \longrightarrow0.
\end{align*}
Define $\zeta_{1},\zeta_{2}$ as in Lemma $\ref{3lemma0}$, we have showed
\begin{align*}
\sup_{|\zeta|\le[\zeta_{1},\zeta_{2}]}\frac{h^{^{\prime\prime\prime}}(\hat
{x}+\zeta)} {h^{^{\prime\prime}}(\hat{x})}\sigma l \longrightarrow0.
\end{align*}
$(\ref{3section201})$ is obtained by using $(\ref{mu3 0060})$. Using
$(\ref{mu3 006})$, for any slowly varying function, it holds
\begin{align*}
h^{^{\prime\prime}}(\hat{x})\sigma^{3} l=\frac{l}{\sqrt{\psi(t)\epsilon(t)t}%
}\longrightarrow0.
\end{align*}
Hence the proof.

\begin{lem}
\label{3lemma02} For $p(x)$ in $(\ref{densityFunction})$, $h\in\mathfrak{R}$,
then for any slowly varying function $l(t)\rightarrow\infty$ as $t\rightarrow
\infty$ such that it holds
\begin{align*}
\sup_{y\in[-l,l]}\frac{|\xi(\sigma y+\hat{x},t)|}{h^{^{\prime\prime}}(\hat
{x})\sigma^{3}}\longrightarrow0,
\end{align*}
where $\xi(x,t)=\epsilon(x,t)+q(x)$.
\end{lem}

Proof: For $y\in[-l,l]$, by $(\ref{3abel01})$ and Lemma $\ref{3lemma01}$ it
holds as $t\rightarrow\infty$
\begin{align}
\label{3lemma0301}\frac{|\epsilon(\sigma y+\hat{x},t)|}{h^{^{\prime\prime}%
}(\hat{x})\sigma^{3}}\le\sup_{|x|\le\sigma l}\Big|\frac{h^{\prime\prime\prime
}(\hat{x}+x)}{h^{\prime\prime}(\hat{x})}\Big|\sigma l\longrightarrow0.
\end{align}
Under condition $(\ref{densityFunction01})$, set $x=\psi(t)$, we get
\begin{align*}
\sup_{|v-\psi(t)|\le\vartheta\psi(t)}|q(v)|\le\frac{1}{\sqrt{t\psi(t)}},
\end{align*}
and it holds for any slowly varying function $l(t)$ as $t\rightarrow\infty$
\begin{align*}
\frac{\sigma l}{\vartheta\psi(t)}=\frac{\sqrt{\psi^{^{\prime}}(t)}
l}{\vartheta\psi(t)}=\sqrt{\frac{\epsilon(t)}{t\psi(t)}}\frac{l}{\vartheta
}\longrightarrow0,
\end{align*}
hence we obtain
\begin{align*}
\sup_{|v-\psi(t)|\le\sigma l}|q(v)|\le\frac{1}{\sqrt{t\psi(t)}}.
\end{align*}
Using this inequality and $(\ref{mu3 006})$, when $y\in[-l,l]$, it holds as
$t\rightarrow\infty$
\begin{align*}
\frac{|q(\sigma y+\hat{x})|}{h^{\prime\prime}(\hat{x})\sigma^{3}}  &
=|q(\sigma y+\hat{x})| \sqrt{\psi(t)\epsilon(t)t}\le\sup_{|v-\psi(t)|\le\sigma
l}|q(v)| \sqrt{\psi(t)\epsilon(t)t}\le\sqrt{\epsilon(t)}\rightarrow0,
\end{align*}
which, together with $(\ref{3lemma0301})$, completes the proof.

\begin{lem}
\label{3lemma1} For $p(x)$ belonging to $(\ref{densityFunction})$,
$h(x)\in\mathfrak{R}$, $\alpha\in\mathbb{N}$, denote by
\begin{align*}
\Psi(t,\alpha):=\int_{0}^{\infty}(x-\hat{x})^{\alpha}e^{tx}p(x)dx,
\end{align*}
then there exists some slowly varying function $l(t)$ such that it holds as
$t\rightarrow\infty$
\begin{align*}
\Psi(t,\alpha)  &  =c\sigma^{\alpha+1} e^{K(\hat{x},t)}T_{1}(t,\alpha
)\big(1+o(1)\big),
\end{align*}
where
\begin{align*}
&  T_{1}(t,\alpha)=\int_{-\frac{l^{1/3}}{\sqrt{2}}}^{\frac{l^{1/3}}{\sqrt{2}}}
y^{\alpha}\exp\big(-\frac{y^{2}}{2}\big)dy-\frac{h^{^{\prime\prime}}(\hat
{x})\sigma^{3}}{6}\int_{-\frac{l^{1/3}}{\sqrt{2}}}^{\frac{l^{1/3}}{\sqrt{2}}}
y^{3+\alpha}\exp\big(-\frac{y^{2}}{2}\big)dy.
\end{align*}

\end{lem}

Proof: By Lemma $\ref{3lemma0}$, for any slowly varying function $l(t)$ it
holds as $t\rightarrow\infty$
\begin{align*}
\sup_{|x-\hat{x}|\le\sigma l} |\epsilon(x,t)|\rightarrow0.
\end{align*}
Given a slowly varying function $l$ with $l(t)\rightarrow\infty$ and define
the interval $I_{t}$ as follows
\begin{align*}
I_{t}:=\Big(-\frac{l^{1/3}\sigma}{\sqrt{2}},\frac{l^{1/3}\sigma}{\sqrt{2}%
}\Big).
\end{align*}
For large enough $\tau$, when $t\rightarrow\infty$ we can partition
$\mathbb{R}_{+}$ as
\begin{align*}
\mathbb{R}_{+}=\{x:0<x< \tau\}\cup\{x:x \in\hat{x}+I_{t}\}\cup\{x:x\ge\tau,x
\notin\hat{x}+I_{t}\},
\end{align*}
where $\tau$ large enough such that it holds for $x>\tau$
\begin{align}
\label{3abeltheorem0010}p(x)< 2ce^{-g(x)}.
\end{align}

Obviously, for fixed $\tau$, $\{x:0<x< \tau\}\cap\{x:x \in\hat{x}+I_{t}\}={Ø}
$ since for large $t$ we have $\min\big(x:x\in\hat{x}+I_{t}\big)\rightarrow
\infty$ as $t\rightarrow\infty$. Hence it holds
\begin{align}
\label{3section1015}\Psi(t,\alpha)  &  =\int_{0}^{\tau}(x-\hat{x})^{\alpha
}e^{tx}p(x)dx+\int_{x \in\hat{x}+I_{t}}(x-\hat{x})^{\alpha}e^{tx}%
p(x)dx+\int_{x \notin\hat{x}+I_{t},x>\tau}(x-\hat{x})^{\alpha}e^{tx}%
p(x)dx\nonumber\\
&  :=\Psi_{1}(t,\alpha)+\Psi_{2}(t,\alpha)+\Psi_{3}(t,\alpha).
\end{align}
We estimate sequentially $\Psi_{1}(t,\alpha),\Psi_{2}(t,\alpha),\Psi
_{3}(t,\alpha)$ in \textbf{Step 1, Step 2 and Step 3}.

\textbf{Step 1:} Using $(\ref{3abeltheorem0010})$, for $\tau$ large enough, we
have
\begin{align}
\label{3section1002}|\Psi_{1}(t,\alpha)|  &  \le\int_{0}^{\tau}|x-\hat
{x}|^{\alpha}e^{tx}p(x)dx\le2c \int_{0}^{\tau}|x-\hat{x}|^{\alpha}%
e^{tx-g(x)}dx\nonumber\\
&  \le2c \int_{0}^{\tau}\hat{x}^{\alpha}e^{tx}dx \le{2c}t^{-1}\hat{x}^{\alpha
}e^{t\tau}.
\end{align}
We show it holds for $h\in\mathfrak{R}$ as $t\rightarrow\infty$
\begin{align}
\label{3section10020}t^{-1}\hat{x}^{\alpha}e^{t\tau}=o(\sigma^{\alpha
+1}e^{K(\hat{x},t)}h^{^{\prime\prime}}(\hat{x})\sigma^{3}).
\end{align}
$(\ref{3section10020})$ is equivalent to
\begin{align*}
\sigma^{-\alpha-4}t^{-1}\hat{x}^{\alpha}e^{t\tau}\big(h^{^{\prime\prime}}%
(\hat{x})\big)^{-1}=o(e^{K(\hat{x},t)}),
\end{align*}
which is implied by
\begin{align*}
\exp\big(-({\alpha}+4)\log\sigma-\log t +\alpha\log\hat{x}+\tau t-\log
h^{^{\prime\prime}}(\hat{x})\big)=o(e^{K(\hat{x},t)}).
\end{align*}
By Lemma $(\ref{3lemma00})$, we know $\log\sigma=o(e^{K(\hat{x},t)})$ as
$t\rightarrow\infty$. So it remains to show $t=o(e^{K(\hat{x},t)})$, $\log
\hat{x}=o(e^{K(\hat{x},t)})$ and $\log h^{^{\prime\prime}}(\hat{x}%
)=o(e^{K(\hat{x},t)})$. Since $\hat{x}=\psi(t)$, it holds
\begin{align}
\label{3section10002}K(\hat{x},t)=t\psi(t)-g(\psi(t))=\int_{1}^{t}
\psi(u)du+\psi(1)-g(1),
\end{align}
where the second equality can be easily verified by the change of variable
$u=h(v)$.

If $h(x)\in R_{\beta}$, by Theorem $(1.5.12)$ of $\cite{Bingham}$, it holds
$\psi(x)\sim x^{1/\beta}l_{1}(x) $ with some slowly varying function
$l_{1}(x)$. $(\ref{3section202})$ and $(\ref{3section10002})$ yield
$t=o(e^{K(\hat{x},t)}) $. In addition, $\log\hat{x}=\log\psi(t)\sim(1/\beta)
\log t=o(e^{K(\hat{x},t)})$. By $(\ref{3lem230})$, it holds $\log
h^{^{\prime\prime}}(\hat{x})=o(t)$. Thus $(\ref{3section10020})$ holds.

If $h(x)\in R_{\infty}$, $\psi(x)\in\widetilde{R_{0}}$ is slowly varying as
$x\rightarrow\infty$. Therefore, by $(\ref{3section203})$ and
$(\ref{3section10002})$, it holds $t=o(e^{K(\hat{x},t)})$ and $\log\hat
{x}=\log\psi(t)=o(e^{K(\hat{x},t)})$. Using $(\ref{mu3 006})$, we have $\log
h^{^{\prime\prime}}(\hat{x})\sim\log t-2\log\hat{x}-2\log\epsilon(t)$. Under
condition $(\ref{3section1030})$, $\log\epsilon(t)=o(t)$, thus it holds $\log
h^{^{\prime\prime}}(\hat{x})=o(t)$. We get $(\ref{3section10020})$.

$(\ref{3section1002})$ and $(\ref{3section10020})$ yield together
\begin{align}
\label{3section1004}|\Psi_{1}(t,\alpha)|=o(\sigma^{\alpha+1}e^{K(\hat{x}%
,t)}h^{^{\prime\prime}}(\hat{x})\sigma^{3}).
\end{align}

\textbf{Step 2:} Notice $\min\big(x:x\in\hat{x}+I_{t}\big)\rightarrow\infty$
as $t\rightarrow\infty$, which implies both $\epsilon(x,t)$ and $q(x)$ go to
$0$ when $x\in\hat{x}+I_{t}$. Using $(\ref{densityFunction})$ and
$(\ref{3abeltheorem001})$, then it holds as $t\rightarrow\infty$
\begin{align*}
\Psi_{2}(t,\alpha)  &  =\int_{x \in\hat{x}+I_{t}}(x-\hat{x})^{\alpha}%
c\exp\big(K(x,t)+q(x)\big)dx\\
&  =\int_{x \in\hat{x}+I_{t}}(x-\hat{x})^{\alpha}c\exp\Big(K(\hat{x}%
,t)-\frac{1}{2}h^{\prime}(\hat{x})\big(x-\hat{x}\big)^{2}\\
&  \qquad\qquad-\frac{1}{6}h^{\prime\prime}(\hat{x})\big(x-\hat{x}%
\big)^{3}+\xi(x,t)\Big)dx,
\end{align*}
where $\xi(x,t)=\epsilon(x,t)+q(x)$. Make the change of variable $y=(x-\hat
{x})/\sigma$, it holds
\begin{align}
\label{3abeltheorem002}\Psi_{2}(t,\alpha)=c\sigma^{\alpha+1} \exp
\big(K(\hat{x},t)\big)\int_{-\frac{l^{1/3}}{\sqrt{2}}}^{\frac{l^{1/3}}%
{\sqrt{2}}} y^{\alpha}\exp\big(-\frac{y^{2}}{2}-\frac{h^{^{\prime\prime}}%
(\hat{x})\sigma^{3}}{6}y^{3}+\xi(\sigma y+\hat{x},t)\big)dy.
\end{align}
On $y\in\big(-{l^{1/3}}/\sqrt{2},{l^{1/3}}/\sqrt{2}\big)$, by
$(\ref{3section2010})$, $|h^{^{\prime\prime}}(\hat{x})\sigma^{3}y^{3}%
|\le|h^{^{\prime\prime}}(\hat{x})\sigma^{3} l|\rightarrow0$ as $t\rightarrow
\infty$. Perform the first order Taylor expansion, it holds as $t\rightarrow
\infty$
\begin{align*}
\exp\big(-\frac{h^{^{\prime\prime}}(\hat{x})\sigma^{3}}{6}y^{3}+\xi(\sigma
y+\hat{x},t)\big)=1-\frac{h^{^{\prime\prime}}(\hat{x})\sigma^{3}}{6}y^{3}%
+\xi(\sigma y+\hat{x},t)+o_{1}(t,y),
\end{align*}
where
\begin{align*}
o_{1}(t,y)=o\big(-\frac{h^{^{\prime\prime}}(\hat{x})\sigma^{3}}{6}y^{3}%
+\xi(\sigma y+\hat{x},t)\big).
\end{align*}
Hence we obtain
\begin{align*}
&  \int_{-\frac{l^{1/3}}{\sqrt{2}}}^{\frac{l^{1/3}}{\sqrt{2}}} y^{\alpha}%
\exp\big(-\frac{y^{2}}{2}-\frac{h^{^{\prime\prime}}(\hat{x})\sigma^{3}}%
{6}y^{3}+\xi(\sigma y+\hat{x},t)\big)dy\\
&  =\int_{-\frac{l^{1/3}}{\sqrt{2}}}^{\frac{l^{1/3}}{\sqrt{2}}} \Big(1-\frac
{h^{^{\prime\prime}}(\hat{x})\sigma^{3}}{6}y^{3}+\xi(\sigma y+\hat{x}%
,t)+o_{1}(t,y)\Big)y^{\alpha}\exp\big(-\frac{y^{2}}{2}\big)dy\\
&  =\int_{-\frac{l^{1/3}}{\sqrt{2}}}^{\frac{l^{1/3}}{\sqrt{2}}} y^{\alpha}%
\exp\big(-\frac{y^{2}}{2}\big)dy-\frac{h^{^{\prime\prime}}(\hat{x})\sigma^{3}%
}{6}\int_{-\frac{l^{1/3}}{\sqrt{2}}}^{\frac{l^{1/3}}{\sqrt{2}}} y^{3+\alpha
}\exp\big(-\frac{y^{2}}{2}\big)dy\\
&  \qquad\qquad\quad+\int_{-\frac{l^{1/3}}{\sqrt{2}}}^{\frac{l^{1/3}}{\sqrt
{2}}} \Big(\xi(\sigma y+\hat{x},t)+o_{1}(t,y)\Big)y^{\alpha}\exp
\big(-\frac{y^{2}}{2}\big)dy.
\end{align*}
Define $T_{1}(t,\alpha)$ and $T_{2}(t,\alpha)$ as follows
\begin{align}
\label{3abeltheorem003} &  T_{1}(t,\alpha)=\int_{-\frac{l^{1/3}}{\sqrt{2}}%
}^{\frac{l^{1/3}}{\sqrt{2}}} y^{\alpha}\exp\big(-\frac{y^{2}}{2}%
\big)dy-\frac{h^{^{\prime\prime}}(\hat{x})\sigma^{3}}{6}\int_{-\frac{l^{1/3}%
}{\sqrt{2}}}^{\frac{l^{1/3}}{\sqrt{2}}} y^{3+\alpha}\exp\big(-\frac{y^{2}}%
{2}\big)dy,\nonumber\\
&  T_{2}(t,\alpha)=\int_{-\frac{l^{1/3}}{\sqrt{2}}}^{\frac{l^{1/3}}{\sqrt{2}}}
\Big(\xi(\sigma y+\hat{x},t)+o_{1}(t,y)\Big)y^{\alpha}\exp\big(-\frac{y^{2}%
}{2}\big)dy.
\end{align}
As for $T_{2}(t,\alpha)$, it holds
\begin{align*}
|T_{2}(t,\alpha)|  &  \le\int_{-\frac{l^{1/3}}{\sqrt{2}}}^{\frac{l^{1/3}%
}{\sqrt{2}}} \Big(|\xi(\sigma y+\hat{x},t)|+|o_{1}(t,y)|\Big)|y|^{\alpha}%
\exp\big(-\frac{y^{2}}{2}\big)dy\\
&  \le\sup_{y\in[-l,l]} |\xi(\sigma y+\hat{x},t)|\int_{-\frac{l^{1/3}}%
{\sqrt{2}}}^{\frac{l^{1/3}}{\sqrt{2}}} |y|^{\alpha}\exp\big(-\frac{y^{2}}%
{2}\big)dy+\int_{-\frac{l^{1/3}}{\sqrt{2}}}^{\frac{l^{1/3}}{\sqrt{2}}}
|o_{1}(t,y)||y|^{\alpha}\exp\big(-\frac{y^{2}}{2}\big)dy\\
&  \le\sup_{y\in[-l,l]} |\xi(\sigma y+\hat{x},t)|\int_{-\frac{l^{1/3}}%
{\sqrt{2}}}^{\frac{l^{1/3}}{\sqrt{2}}} |y|^{\alpha}\exp\big(-\frac{y^{2}}%
{2}\big)dy\\
&  \qquad\qquad+\int_{-\frac{l^{1/3}}{\sqrt{2}}}^{\frac{l^{1/3}}{\sqrt{2}}}
\Big(\big|o\big(\frac{h^{^{\prime\prime}}(\hat{x})\sigma^{3}}{6}%
y^{3}\big)\big|+\big|o\big(\xi(\sigma y+\hat{x},t)\big)\big|\Big)|y|^{\alpha
}\exp\big(-\frac{y^{2}}{2}\big)dy\\
&  \le2\sup_{y\in[-l,l]} |\xi(\sigma y+\hat{x},t)|\int_{-\frac{l^{1/3}}%
{\sqrt{2}}}^{\frac{l^{1/3}}{\sqrt{2}}} |y|^{\alpha}\exp\big(-\frac{y^{2}}%
{2}\big)dy+|o(h^{^{\prime\prime}}(\hat{x})\sigma^{3})|\int_{-\frac{l^{1/3}%
}{\sqrt{2}}}^{\frac{l^{1/3}}{\sqrt{2}}} |y|^{3+\alpha} \exp\big(-\frac{y^{2}%
}{2}\big)dy\\
&  = |o(h^{^{\prime\prime}}(\hat{x})\sigma^{3})|\Big(\int_{-\frac{l^{1/3}%
}{\sqrt{2}}}^{\frac{l^{1/3}}{\sqrt{2}}} |y|^{\alpha}\exp\big(-\frac{y^{2}}%
{2}\big)dy+\int_{-\frac{l^{1/3}}{\sqrt{2}}}^{\frac{l^{1/3}}{\sqrt{2}}}
|y|^{3+\alpha} \exp\big(-\frac{y^{2}}{2}\big)dy\Big),
\end{align*}
where last equality holds from Lemma $\ref{3lemma02}$. Since the integrals in
the last equality are both bounded, it holds as $t\rightarrow\infty$
\begin{align*}
T_{2}(t,\alpha)=o(h^{^{\prime\prime}}(\hat{x})\sigma^{3}).
\end{align*}

When $\alpha$ is even, the second term of $T_{1}(t,\alpha)$ vanishes. When
$\alpha$ is odd, the first term of $T_{1}(t,\alpha)$ vanishes. Obviously,
$T_{1}(t,\alpha)$ is at least the same order than $h^{^{\prime\prime}}(\hat
{x})\sigma^{3}$. Therefore it follows as $t\rightarrow\infty$
\begin{align}
\label{3abeltheorem004}T_{2}(t,\alpha)=o(T_{1}(t,\alpha)).
\end{align}
Using $(\ref{3abeltheorem002})$, $( \ref{3abeltheorem003})$ and
$(\ref{3abeltheorem004})$ we get
\begin{align}
\label{3section1014}\Psi_{2}(t,\alpha)  &  =c\sigma^{\alpha+1} \exp
\big(K(\hat{x},t)\big)T_{1}(t,\alpha)\big(1+o(1)\big).
\end{align}

\textbf{Step 3:} Given $h\in\mathfrak{R}$, for any $t$, $K(x,t)$ as a function
of $x$ ($x>\tau$) is concave since
\begin{align*}
K^{\prime\prime}(x,t)=-h^{\prime}(x)<0.
\end{align*}
Thus we get for $x\notin\hat{x}+I_{t}$ and $x>\tau$
\begin{align}
\label{3abeltheorem005}K(x,t)-K(\hat{x},t)\le\frac{K(\hat{x}+\frac
{l^{1/3}\sigma}{\sqrt{2}}sgn(x-\hat{x}),t)-K(\hat{x},t)}{\frac{l^{1/3}\sigma
}{\sqrt{2}}sgn(x-\hat{x})}(x-\hat{x}),
\end{align}
where
\begin{align*}
sgn(x-\hat{x})=
\begin{cases}
1 \qquad\;\;\; if \quad x\ge\hat{x},\\
-1 \qquad if \quad x<\hat{x}.
\end{cases}
\end{align*}
Using $(\ref{3abeltheorem001})$, we get
\begin{align*}
{K(\hat{x}+\frac{l^{1/3}\sigma}{\sqrt{2}}sgn(x-\hat{x}),t)-K(\hat{x},t)}
\le-\frac{1}{8}h^{\prime}(\hat{x})l^{2/3}\sigma^{2}=-\frac{1}{8}l^{2/3},
\end{align*}
which, combined with $(\ref{3abeltheorem005})$, yields
\begin{align*}
K(x,t)-K(\hat{x},t)\le-\frac{\sqrt{2}}{8}l^{1/3}\sigma^{-1}|x-\hat{x}|.
\end{align*}
\newline We obtain
\begin{align*}
|\Psi_{3}(t,\alpha)|  &  \le2c\int_{x \notin\hat{x}+I_{t},x>\tau}|x-\hat
{x}|^{\alpha}\exp\big(K(x,t)\big)dx\\
&  \le2c\int_{|x-\hat{x}|>\frac{l^{1/3}\sigma}{\sqrt{2}}}|x-\hat{x}|^{\alpha
}\exp\big (K(x,t)\big)dx\\
&  \le2 ce^{K(\hat{x},t)}\int_{|x- \hat{x}|>\frac{l^{1/3}\sigma}{\sqrt{2}}%
}|x-\hat{x}|^{\alpha}\exp\big(-\frac{\sqrt{2}}{8}l^{1/3}\sigma^{-1}|x-\hat
{x}|\big)dx\\
&  =2c e^{K(\hat{x},t)}\sigma^{\alpha+1}\int_{|y|>\frac{l^{1/3}}{\sqrt{2}}%
}|y|^{\alpha}\exp\big(-\frac{\sqrt{2}}{8}l^{1/3}|y|\big)dy\\
&  =2c e^{K(\hat{x},t)}\sigma^{\alpha+1}\int_{|y|>\frac{l^{1/3}}{\sqrt{2}}}
\exp\big(-\frac{\sqrt{2}}{8}l^{1/3}|y|+\alpha\log|y|\big)dy\\
&  =2c e^{K(\hat{x},t)}\sigma^{\alpha+1} \Big (2e^{-l^{2/3}/8}%
\big(1+o(1)\big)\Big),
\end{align*}
where last equality holds when $l\rightarrow\infty$ (see e.g. Theorem 4.12.10
of \cite{Bingham}). With $(\ref{3section1014})$, we obtain
\begin{align*}
\Big|\frac{\Psi_{3}(t,\alpha)}{\Psi_{2}(t,\alpha)}\Big|\le\frac{8e^{-l^{2/3}%
/8}}{|T_{1}(t,\alpha)|}.
\end{align*}
In \textbf{Step 2}, we know $T_{1}(t,\alpha)$ has at least the order
$h^{^{\prime\prime}}(\hat{x})\sigma^{3}$. Hence there exists some positive
constant $Q$ and $l_{2}(t)\rightarrow\infty$ such that it holds as
$t\rightarrow\infty$
\begin{align*}
\Big|\frac{\Psi_{3}(t,\alpha)}{\Psi_{2}(t,\alpha)}\Big|  &  \le\frac
{Qe^{-l_{2}^{2/3}/8}}{h^{^{\prime\prime}}(\hat{x})\sigma^{3}}.
\end{align*}
For example, we can take $l_{2}(t)=(\log t)^{3}$.

If $h\in R_{\beta}$, by $(\ref{3hfe1})$, it is easy to know $h^{^{\prime
\prime}}(\hat{x})\sigma^{3} \ge{1}/{t^{1+1/(2\beta)}}$, thus we have
\begin{align*}
\Big|\frac{\Psi_{3}(t,\alpha)}{\Psi_{2}(t,\alpha)}\Big|  &  \le Q\exp
\big(-l_{2}^{2/3}/8+(1+1/(2\beta))\log t\big)\longrightarrow0.
\end{align*}

If $h\in R_{\infty}$, using $(\ref{mu3 006})$, then it holds as $t\rightarrow
\infty$
\begin{align}
\label{3section1016}\Big|\frac{\Psi_{3}(t,\alpha)}{\Psi_{2}(t,\alpha)}\Big|
&  \le2Q\exp\big(-l_{2}^{2/3}/8+\log\sqrt{t\psi(t)\epsilon(t)}\big)\nonumber\\
&  =2Q\exp\Big(-l_{2}^{2/3}/8+({1}/{2})\big(\log t+\log\psi(t)+\log
\epsilon(t)\big)\Big)\nonumber\\
&  \longrightarrow0,
\end{align}
where last line holds since $\log\psi(t)=O(\log t)$. The proof is completed by
combining $(\ref{3section1015})$, $(\ref{3section1004})$, $(\ref{3section1014}%
)$ and $(\ref{3section1016})$.

\bigskip\textbf{Proof of Theorem $\ref{order of s}$:} By Lemma $\ref{3lemma1}%
$, if $\alpha=0$, it holds $T_{1}(t,0)\sim\sqrt{2\pi}$ as $t\rightarrow\infty
$, hence for $p(x)$ defined in $(\ref{densityFunction})$, we can approximate
$X$'s moment generating function $\Phi(t)$
\begin{align}
\label{3moment001}\Phi(t)=\int_{0}^{\infty}e^{tx}p(x)dx=c\sqrt{2\pi}\sigma
e^{K(\hat{x},t)}\big(1+o(1)\big).
\end{align}

If $\alpha=1$, it holds as $t\rightarrow\infty$,
\begin{align*}
&  T_{1}(t,1)=-\frac{h^{^{\prime\prime}}(\hat{x})\sigma^{3}}{6}\int
_{-\frac{l^{1/3}}{\sqrt{2}}}^{\frac{l^{1/3}}{\sqrt{2}}} y^{4}\exp
\big(-\frac{y^{2}}{2}\big)dy=-\frac{\sqrt{2\pi}h^{^{\prime\prime}}(\hat
{x})\sigma^{3}}{2}\big(1+o(1)\big),
\end{align*}
hence we have with $\Psi(t,\alpha)$ defined in Lemma $\ref{3lemma1}$
\begin{align}
\label{3moment020}\Psi(t,1)  &  =-c\sqrt{2\pi}\sigma^{2} e^{K(\hat{x},t)}%
\frac{h^{^{\prime\prime}}(\hat{x})\sigma^{3}}{2}\big(1+o(1)\big)=-\Phi
(t)\frac{h^{^{\prime\prime}}(\hat{x})\sigma^{4}}{2}\big(1+o(1)\big),
\end{align}
which, together with the definition of $\Psi(t,\alpha)$, yields
\begin{align}
\label{3moment02}\int_{0}^{\infty}xe^{tx}p(x)dx=\Psi(t,1)+\hat{x}%
\Phi(t)=\Big(\hat{x}-\frac{h^{^{\prime\prime}}(\hat{x})\sigma^{4}}%
{2}\big(1+o(1)\big)\Big)\Phi(t).
\end{align}
Hence we get
\begin{align}
\label{3moment021}m(t)=\frac{d\log\Phi(t)}{dt}  &  =\frac{\int_{0}^{\infty
}xe^{tx}p(x)dx}{\Phi(t)}=\hat{x}-\frac{h^{^{\prime\prime}}(\hat{x})\sigma^{4}%
}{2}\big(1+o(1)\big).
\end{align}

Set $\alpha=2$, as $t\rightarrow\infty$, it follows
\begin{align}
\label{3moment03}\Psi(t,2)  &  =c\sigma^{3} e^{K(\hat{x},t)}\int
_{-\frac{l^{1/3}}{\sqrt{2}}}^{\frac{l^{1/3}}{\sqrt{2}}} y^{2} \exp
\big(-\frac{y^{2}}{2}\big)dy\big(1+o(1)\big)\nonumber\\
&  =c\sqrt{2\pi}\sigma^{3} e^{K(\hat{x},t)}\big(1+o(1)\big)=\sigma^{2}%
\Phi(t)\big(1+o(1)\big).
\end{align}
Using $(\ref{3moment020}),(\ref{3moment021})$ and $(\ref{3moment03})$, we
have
\begin{align*}
&  \int_{0}^{\infty}\big(x-m(t)\big)^{2} e^{tx}p(x)dx =\int_{0}^{\infty
}\big(x-\hat{x}+\hat{x}-m(t)\big)^{2} e^{tx}p(x)dx\\
&  =\int_{0}^{\infty}\big(x-\hat{x}\big)^{2} e^{tx}p(x)dx +2\big(\hat
{x}-m(t)\big)\int_{0}^{\infty}(x-\hat{x}) e^{tx}p(x)dx+\big(\hat
{x}-m(t)\big)^{2}\Phi(t)\\
&  =\Psi(t,2)+2\big(\hat{x}-m(t)\big)\Psi(t,1)+\big(\hat{x}-m(t)\big)^{2}%
\Phi(t)\\
&  =\sigma^{2}\Phi(t)\big(1+o(1)\big)-h^{^{\prime\prime}}(\hat{x})\sigma
^{4}\Big(\Phi(t)\frac{h^{^{\prime\prime}}(\hat{x})\sigma^{4}}{2}%
\Big)\big(1+o(1)\big) +\Big(\frac{h^{^{\prime\prime}}(\hat{x})\sigma^{4}}%
{2}\Big)^{2}\Phi(t)\big(1+o(1)\big)\\
&  =\Big(\sigma^{2}-\frac{(h^{^{\prime\prime}}(\hat{x})\sigma^{4})^{2}}%
{4}\Big)\Phi(t)\big(1+o(1)\big),
\end{align*}
thus we have
\begin{align}
\label{3moment02012}s^{2}(t)=\frac{d^{2}\log\Phi(t)}{dt^{2}}  &  =\frac
{\int_{0}^{\infty}\big(x-m(t)\big)^{2}e^{tx}p(x)dx}{\Phi(t)}=\Big(\sigma
^{2}-\frac{(h^{^{\prime\prime}}(\hat{x})\sigma^{4})^{2}}{4}%
\Big)\big(1+o(1)\big).
\end{align}

Set $\alpha=3$, the first term of $T_{1}(t,3)$ vanishes, we obtain as
$t\rightarrow\infty$
\begin{align}
\label{3moment0202}\Psi(t,3)  &  =-c\sqrt{2\pi}\sigma^{4} e^{K(\hat{x}%
,t)}\frac{h^{^{\prime\prime}}(\hat{x})\sigma^{3}}{2}\int_{-\frac{l^{1/3}%
}{\sqrt{2}}}^{\frac{l^{1/3}}{\sqrt{2}}} \frac{1}{\sqrt{2\pi}}y^{6}%
\exp\big(-\frac{y^{2}}{2}\big)dy\nonumber\\
&  =-cM_{6}\sqrt{2\pi}e^{K(\hat{x},t)}\frac{h^{^{\prime\prime}}(\hat{x}%
)\sigma^{7}}{2}\big(1+o(1)\big)=-M_{6}\frac{h^{^{\prime\prime}}(\hat{x}%
)\sigma^{6}}{2}\Phi(t)\big(1+o(1)\big),
\end{align}
where $M_{6}$ denotes the sixth order moment of standard normal distribution.
Using $(\ref{3moment020}),(\ref{3moment021})$, $(\ref{3moment03})$ and
$(\ref{3moment0202})$, we have as $t\rightarrow\infty$
\begin{align*}
&  \int_{0}^{\infty}\big(x-m(t)\big)^{3} e^{tx}p(x)dx =\int_{0}^{\infty
}\big(x-\hat{x}+\hat{x}-m(t)\big)^{3} e^{tx}p(x)dx\\
&  =\int_{0}^{\infty}\Big((x-\hat{x})^{3}+3(x-\hat{x})^{2}\big(\hat
{x}-m(t)\big)+3(x-\hat{x})\big(\hat{x}-m(t)\big)^{2} +\big(\hat{x}%
-m(t)\big)^{3}\Big)e^{tx}p(x)dx\\
&  =\Psi(t,3)+3\big(\hat{x}-m(t)\big)\Psi(t,2)+3\big(\hat{x}-m(t)\big)^{2}%
\Psi(t,1)+\big(\hat{x}-m(t)\big)^{3}\Phi(t)\\
&  =-M_{6}\frac{h^{^{\prime\prime}}(\hat{x})\sigma^{6}}{2}\Phi
(t)\big(1+o(1)\big)+(3/2)h^{^{\prime\prime}}(\hat{x})\sigma^{4}(\sigma^{2}%
\Phi(t))\big(1+o(1)\big)\\
&  \qquad-3\Big(\frac{h^{^{\prime\prime}}(\hat{x})\sigma^{4}}{2}\Big)^{2}%
\Phi(t)\frac{h^{^{\prime\prime}}(\hat{x})\sigma^{4}}{2}%
\big(1+o(1)\big)+\Big(\frac{h^{^{\prime\prime}}(\hat{x})\sigma^{4}}%
{2}\Big)^{3}\Phi(t)\big(1+o(1)\big)\\
&  =\Big(\frac{3-M_{6}}{2}h^{^{\prime\prime}}(\hat{x})\sigma^{6}%
-\frac{(h^{^{\prime\prime}}(\hat{x})\sigma^{4})^{3}}{4}\Big)\Phi
(t)\big(1+o(1)\big),
\end{align*}
hence we get
\begin{align}
\label{3moment02022}\mu_{3}(t)=\frac{d^{3}\log\Phi(t)}{dt^{3}}  &  =\frac
{\int_{0}^{\infty}\big(x-m(t)\big)^{3}e^{tx}p(x)dx}{\Phi(t)}=\Big(\frac
{3-M_{6}}{2}h^{^{\prime\prime}}(\hat{x})\sigma^{6}-\frac{(h^{^{\prime\prime}%
}(\hat{x})\sigma^{4})^{3}}{4}\Big)\big(1+o(1)\big).
\end{align}

Finally, we finish the proof by simplifying $(\ref{3moment021})$
$(\ref{3moment02012})$ and $(\ref{3moment02022})$.

\textbf{Case 1:} $h\in R_{\beta}$. We have gotten in $(\ref{3lem230})$
\begin{align*}
h^{^{\prime\prime}}(\hat{x})=\big(\beta(\beta-1)+o(1)\big)\psi(t)^{\beta-2}
l_{0}(\psi(t)),
\end{align*}
where $l_{0} \in R_{0}$. In $(\ref{mu3 1013})$, we have $\sigma^{2}\sim
\psi(t)/(\beta t)$, hence it holds
\begin{align*}
h^{\prime\prime}(\hat{x})\sigma^{4}=\frac{\beta-1}{\beta}\frac{\psi(t)^{\beta
}}{t^{2}}l_{0}(\psi(t))\big(1+o(1)\big)=\frac{\beta-1}{\beta}\frac{l_{0}%
(\psi(t))l_{1}(t)^{\beta}}{t}\big(1+o(1)\big),
\end{align*}
where last equality holds since $\psi(t)\sim t^{1/\beta}l_{1}(t)$ for some
slowly varying function $l_{1}$. Obviously, $h^{\prime\prime}(\hat{x}%
)\sigma^{4}=o(\hat{x})$, thus we have
\begin{align*}
m(t)\sim\hat{x}=\psi(t).
\end{align*}
It holds also as $t\rightarrow\infty$
\begin{align*}
\frac{(h^{\prime\prime}(\hat{x})\sigma^{4})^{2}}{\sigma^{2}}={(\beta-1)^{2}%
}\frac{l_{0}(\psi(t))^{2}}{\psi(t)^{2}}\big(1+o(1)\big)\longrightarrow0,
\end{align*}
which implies $(h^{\prime\prime}(\hat{x})\sigma^{4})^{2}=o(\sigma^{2})$.
Therefore it follows
\begin{align}
\label{3appen01}s^{2}(t)\sim\sigma^{2}=\psi^{^{\prime}}(t).
\end{align}
For $\mu_{3}$, it holds $(h^{^{\prime\prime}}(\hat{x})\sigma^{4}%
)^{3}=o(h^{^{\prime\prime}}(\hat{x})\sigma^{6})$ since
\begin{align*}
\frac{(h^{^{\prime\prime}}(\hat{x})\sigma^{4})^{3}}{h^{^{\prime\prime}}%
(\hat{x})\sigma^{6}}=h^{^{\prime\prime}}(\hat{x})^{2}\sigma^{6} =\frac
{(\beta-1)^{2}}{\beta}\frac{\psi(t)^{2\beta-1}l_{0}(\psi(t))^{2}}{t^{3}%
}\big(1+o(1)\big)\longrightarrow0,
\end{align*}
where last step holds from the fact $\psi(t)^{2\beta-1}/t^{3} \sim
l_{1}(t)^{2\beta-1}/t^{1+1/\beta}$. We have
\begin{align}
\label{3appen02}\mu_{3}(t)\sim\frac{3-M_{6}}{2}h^{^{\prime\prime}}(\hat
{x})\sigma^{6}.
\end{align}
It is straightforward that $(\ref{mu3 0010})$ holds for $h\in R_{\beta}$, thus
$h^{^{\prime\prime}}(\hat{x})\sigma^{6}=-\psi^{^{\prime\prime}}(t)/(\psi
^{^{\prime}}(t))^{3}*(\psi^{^{\prime}}(t))^{3}=-\psi^{^{\prime\prime}}(t)$ .
We get
\begin{align*}
\mu_{3}(t)\sim\frac{M_{6}-3}{2}\psi^{^{\prime\prime}}(t).
\end{align*}

\textbf{Case 2:} If $h\in R_{\infty}$, recall that we have obtained in
$(\ref{mu3 006})$
\begin{align*}
h^{\prime\prime}(\hat{x})=-\frac{\psi^{\prime\prime}(t)}{\big(\psi^{\prime
}(t)\big)^{3}} =\frac{t}{\psi^{2}(t)\epsilon^{2}(t)}\big(1+o(1)\big),
\end{align*}
consider $\sigma^{2}=\psi^{^{\prime}}(t)=\psi(t)\epsilon(t)/t$, it holds
\begin{align*}
h^{\prime\prime}(\hat{x})\sigma^{4}=\frac{1}{t}\big(1+o(1)\big).
\end{align*}
Notice $h^{\prime\prime}(\hat{x})\sigma^{4}=o(\hat{x})$ as $t\rightarrow
\infty$, hence it holds
\begin{align*}
m(t)\sim\hat{x}=\psi(t).
\end{align*}
And as $t\rightarrow\infty$ it holds $\big(h^{\prime\prime}(\hat{x})\sigma
^{4}\big)^{2}\sim1/t^{2}=o(\sigma^{2})$, thus we obtain
\begin{align*}
s^{2}(t)\sim\sigma^{2}=\psi^{^{\prime}}(t).
\end{align*}
As regards to $\mu_{3}(t)$, we have $\big(h^{\prime\prime}(\hat{x})\sigma
^{4}\big)^{3}\sim1/t^{3}$, but $h^{\prime\prime}(\hat{x})\sigma^{6} \sim
\psi(t)\epsilon(t)/t^{2}$, hence it holds $\big(h^{\prime\prime}(\hat
{x})\sigma^{4}\big)^{3}=o\big(h^{\prime\prime}(\hat{x})\sigma^{6}\big)$. It
follows
\begin{align*}
\mu_{3}(t)\sim\frac{M_{6}-3}{2}\psi^{^{\prime\prime}}(t).
\end{align*}
\newline\textbf{Proof of Corollary $\ref{3cor1}:$} \textbf{Case 1:} If $h\in
R_{\beta}$. By $(\ref{3appen01})$ and $(\ref{3appen02})$, it holds as
$t\rightarrow\infty$
\begin{align}
\label{3appen03}\frac{\mu_{3}}{s^{3}}\sim\frac{M_{6}-3}{2}h^{^{\prime\prime}%
}(\hat{x})\sigma^{3}.
\end{align}
Then using $(\ref{3lem230})$ and $(\ref{mu3 1013})$, we get for $l_{0}\in
R_{0}$
\begin{align}
\label{3appen04}h^{^{\prime\prime}}(\hat{x})\sigma^{3}  &  \sim\beta
(\beta-1)\psi(t)^{\beta-2} l_{0}(\psi(t))\Big(\frac{\psi(t)}{\beta
t}\Big)^{3/2}\nonumber\\
&  =\frac{\beta-1}{\sqrt{\beta}} l_{0}(\psi(t))\frac{\psi(t)^{\beta-1/2}}{
t^{3/2}}\longrightarrow0,
\end{align}
where last step holds since $\psi(t)\sim t^{1/\beta}l_{1}(t)$ for some slowly
varying function $l_{1}(t)$. $(\ref{3appen03})$ and $(\ref{3appen04})$ yields
$(\ref{3cor1})$.

\textbf{Case 2:} If $h\in R_{\infty}$. In \textbf{(1)} we have showed it
holds
\begin{align*}
\frac{\mu_{3}(t)}{s^{3}(t)}\sim\frac{M_{6}-3}{2}\frac{\psi^{^{\prime\prime}%
}(t)}{\psi^{\prime}(t)^{3/2}}.
\end{align*}
By $(\ref{mu3 002})$ and $(\ref{mu3 003})$, we have as $t\rightarrow\infty$
\begin{align*}
\frac{\psi^{^{\prime\prime}}(t)}{\psi^{\prime}(t)^{3/2}}\sim-\frac
{\psi(t)\epsilon(t)}{t^{2}}\Big(\frac{\psi(t)\epsilon(t)}{t}\Big)^{-3/2}
=-\frac{1}{\sqrt{t\psi(t)\epsilon(t)}}\longrightarrow0,
\end{align*}
where last step holds under condition $(\ref{3section1030})$. Hence the proof.

\end{document}